\documentclass[11pt]{article}
\usepackage{graphics, epic,eepic, url}

\newcounter{prob}
\newenvironment{problem}{\begin{list}{{\it Problem \arabic{prob}
        :}}{\usecounter{prob}
      \settowidth{\labelwidth}{\it Problem 5 :}
      \setlength{\labelsep}{0.2cm}
      \setlength{\leftmargin}{\labelwidth}
      \addtolength{\leftmargin}{\labelsep}}}{\end{list}}
        
\newcommand{\twoo}{\mbox{\rm $2$--orb}}
\newcommand{\aut}{\mbox{\rm Aut}}
\newcommand{\COCO}{{\sf COCO}}

\def\0{\,(0)\,}
\def\1{\,(1)\,}

\begin{document}
\begin{titlepage}

\begin{center}
  {\LARGE Algebraic Combinatorics in Mathematical Chemistry.\\
   Methods and Algorithms. \\
  II. Program Implementation of the Weisfeiler-Leman Algorithm\\
\vspace{10mm}
  {\normalsize \url{arXiV.org} {\em Version}}}\\
\vspace{20mm}

  {\Large Luitpold Babel$^{\mbox{\small 1}}$}\\ 
  {\sl\normalsize Institut f\"ur Mathematik}\\ 
  {\sl\normalsize Technische Universit\"at M\"unchen}\\
  {\sl\normalsize D--80290 M\"unchen, Germany}\\
\vspace{2mm}
  {\Large Irina V. Chuvaeva$^{\mbox{\small 2}}$}\\ 
  {\sl\normalsize Laboratory of Mathematical Chemistry}\\ 
  {\sl\normalsize N.D.Zelinski\v{\i} Institute of Organic Chemistry}\\
  {\sl\normalsize Moscow, Russia}\\
\vspace{2mm}
  {\Large Mikhail Klin$^{\mbox{\small 1}}$}\\ 
  {\sl\normalsize Department of Mathematics}\\ 
  {\sl\normalsize Ben-Gurion University of the Negev}\\
  {\sl\normalsize 84105 Beer-Sheva, Israel}\\
\vspace{2mm}
  {\Large Dmitrii V. Pasechnik$^{\mbox{\small 3}}$}\\ 
  {\sl\normalsize Laboratory of Discrete Mathematics, VNIISI}\\
  {\sl\normalsize Moscow, Russia}\\
\end{center}

\vspace*{45mm}
{\footnotesize
\hspace{3mm}  $^{\mbox{\tiny 1}}$Supported by the grant
No. I-0333-263.06/93 from the G.I.F., the German-Israeli
Foundation for Scientific Research and Development \\
\hspace*{3mm}  $^{\mbox{\tiny 2}}$Current address: 249073, selo Nedel'noe,
Maloyaroslavetski\v{\i} raion, Kalu\v{z}skaya oblast', Russia\\
\hspace*{3mm}  $^{\mbox{\tiny 3}}$Current affiliation: Division of 
Mathematical Sciences, SPMS, Nanyang Technological University, 21
Nanyang Link, 637371 Singapore. URL: 
\url{http://www.ntu.edu.sg/home/dima/}
}

\end{titlepage}
\begin{abstract}
The stabilization algorithm of Weisfeiler and Leman has as an input any
square matrix $A$ of order $n$ and returns the minimal cellular
(coherent) algebra $W(A)$ which includes $A$.

In case when $A=A(\Gamma)$ is the adjacency matrix of a graph $\Gamma$
the algorithm examines all configurations in $\Gamma$ having three
vertices and, according to this information, partitions vertices and
ordered pairs of vertices into equivalence classes. The resulting
construction allows to associate to each graph $\Gamma$ a matrix
algebra $W(\Gamma):= W\left(A(\Gamma)\right)$ which is an invariant
of the graph $\Gamma$. For many classes of graphs, in particular for
most of the molecular graphs, the algebra $W(\Gamma)$ coincides with
the centralizer algebra of the automorphism group $\aut(\Gamma)$. In
such a case the partition returned by the stabilization algorithm 
is equal to the partition into orbits of $\aut(\Gamma)$.

We give algebraic and combinatorial descriptions of the
Weisfeiler--Leman algorithm and present an efficient computer
implementation of the algorithm written in C. The results obtained by
testing the program on a considerable number of examples of
graphs, in particular on some chemical molecular graphs, are also
included.
\end{abstract}

\section{Introduction}
The fundamental problem of graph symmetry perception arises in
numerous areas of chemistry and physics. In this context molecules
are modelled by graphs where the vertices represent atoms and the
edges represent bonds.  The aim is to find equivalence classes of
``elements'' (e.g.\ vertices, edges, pairs of vertices, subgraphs etc.) of
the graph, or more rigorously, the orbits of the action of the
automorphism group $\aut(\Gamma)$ of the graph $\Gamma$ on the set of all
``elements'' of a considered ``mode'' (for example the orbits of $\aut(\Gamma)$
on the set of vertices). It is clear that in order to use such a statement
one has to get at some intermediate stage a convenient description of $\aut(\Gamma)$.

However, usually chemists avoid the computation of $\aut(\Gamma)$,
preferring to use certain invariants in order to get a classification
of ``elements'' of the graph.

Let us express the above claim more accurately for the case of vertices.
A vertex invariant in a graph $\Gamma$ is a property or a parameter of
vertices which is preserved by any of its automorphisms, i.e.
the property does not depend on the
labelling of the graph.  Let $\phi$ be a function which is defined on
the set of all vertices of $\Gamma$. Then $\phi$ \ is an invariant of vertices if
$\phi(x)=\phi(x^\prime)$ whenever the vertices $x$, $x^\prime$ belong
to the same orbit of $\aut(\Gamma)$.

Let us now consider a certain invariant (a set of invariants). Then we may define
that two vertices belong to the same equivalency class if and only if they
have the same value of the invariant (invariants) being taken into account.
Whatever the classification approach, the resulting
partition cannot be finer than the partition into the orbits of the
automorphism group (the automorphism partition).

From an algorithmic point of view the main goal in classification of
vertices is to find an algorithm which ensures to produce the automorphism
partition of a graph and which is known to be theoretically
efficient (this means that the running time is restricted by a
polynomial in the size of the input). However, all known polynomial time methods for
the automorphism partitioning problem yield heuristic solutions, i.e.\ 
they result in partitions where the equivalence classes are orbits or
unions of orbits. 

It is well known (see e.g.\ \cite{ReaC77},\ \cite{Pon94c}) that the automorphism
partitioning problem (which is strongly related to the problem of
graph symmetry\footnote[1]{ In this paper we are not interested in
  finding the automorphism group of a graph (what is commonly meant by
  solving the problem of graph symmetry) but rather in finding the
  automorphism partition.})  is algorithmically equivalent to the
graph isomorphism problem (the problem of graph identification) in the
following sense.  Assume that a polynomial algorithm is given solving
one of the two problems.  Then it is possible to construct from this
algorithm a second polynomial algorithm which solves the other
problem. However, both problems seem to be hard from a 
computational point of view. Therefore heuristic
approaches are still considered as practically helpful.

The simplest of these approaches are restricted to a classification of
vertices (in chemical terms they determine atom equivalence only) and
are based on different techniques, often applied iteratively, using
the valencies of the vertices.  The major weakness of these methods is
that they will not give any partition for regular graphs (graphs where
all vertices have the same valencies).  For that reason it seems quite
natural to extend these techniques to a classification of vertices and
edges (i.e.\ atoms and pairs of atoms).  Now, not only configurations
of two vertices (which define the valencies) but also configurations
consisting of three vertices have to be considered.  This is exactly
the basic idea of the algorithm of Weisfeiler and Leman.  It
partitions all vertices and ordered pairs of vertices of a graph by
examining all ordered triples of vertices. This approach turns out to
be the most powerful method in some class of graph symmetry perception
algorithms.

The algebraic object which is constructed in this way and which has been introduced 
by Weis\-feiler and Leman in \cite{WeiL68} is a cellular algebra.
It is an invariant of the underlying graph. 
Under a different point of view (without any relation to the automorphism partitioning 
problem) this object has also been found and called coherent configuration by Higman 
in \cite{Hig70}. In greater detail, Weisfeiler-Leman's approach was described (in 
English) in \cite{Wei76}. However, during approximately twenty years, their ideas 
were completely unknown and neglected in mathematical chemistry. Nowadays, the 
approach itself and its interrelations with the identification and symmetry 
perception of graphs are rather familiar to the experts in algebraic combinatorics
(see e.g.\ \cite{Hig87}, \cite{Fri89}, \cite{FarIK90}, \cite{Pon93b},
\cite{Pon94b}, \cite{FarKM94}),
however a careful investigation of its abilities still remains a topical problem.

More or less the same as Weisfeiler-Leman's approach was independently elaborated
by G.\ Tinhofer (partly in joint work with his student J.~Hinteregger)
in \cite{Tin75} and \cite{HinT77}, however 
without explicitly describing the resulting algebraic object.
In 1989, Ch.\ and G.\ R\"ucker (a chemist and a mathematician) realized
the necessity of having a method for the partition of atom pairs, and produced
a heuristic computer program for that purpose by simple reasoning without using
any group theoretical machinery \cite{RueR90a}, \cite{RueR90b}, \cite{RueR91}.

The main purpose of our paper is not only to draw attention to the
algorithm of Weisfeiler and Leman but, mainly, to present a good and
practical program implementation. To our knowledge, no such
implementation exists up to now.  Only few attempts were made in the
past. A first version of our program has been described by I.V.\ 
Chuvaeva, M.\ Klin and D.V.\ Pasechnik in \cite{ChuKP92}.  By means of
a very careful revision we now realized all advantages and
disadvantages of this implementation.  Recently, L.\ Babel established
in \cite{Bab95} the theoretical complexity of the Weisfeiler-Leman
algorithm and, using these considerations, also created a computer
program. A detailed description is given in \cite{BabBLT97}.  Our
common experiences (including some important suggestions of Ch. Pech)
enabled us to modify the version of the program by Chuvaeva et al.\ into a very
fast program implementation. A comparison with the program of Babel
et al.\ shows that our program, although inferior with respect to
theoretical complexity, is much more efficient from a practical point
of view.

This paper is organized as follows. In Section 2 we
introduce the standard terminology and some basic definitions.  After
that a brief survey on previous approaches to graph stabilization is
given in Section 3.  Section 4 contains the definitions of a cellular
algebra and related algebraic objects, states some important
properties and interpretations and introduces the cellular algebra
which is associated to a given graph.  In Section 5 an algebraic
description of the algorithm of Weisfeiler and Leman is presented,
followed in Section 6 by a more illustrating combinatorial
interpretation.  Thus, the algorithm is exposed from two different
points of view, the first using matrix notation, the second using
graph theoretical notation.  Sections 7 and 8 give a description of
the program implementation and an estimation of its complexity.
Furthermore, in a brief excursion we present the main ideas of the
complexity considerations and of the implementation of Babel's
algorithm.  Finally, in Section 9, extended testing of our program on
a large number of examples is documented in order to demonstrate its
capability.  We conclude with a discussion in Section 10.

This work is the second contribution in a series of papers
\cite{KliRRT95}, \cite{FurKT} concerning different aspects of
algebraic combinatorics with emphasis on applications in mathematical
chemistry.  The series introduces the basic concepts of algebraic
combinatorics and presents some of the main features and tools for
perception of symmetry properties of combinatorial objects.  Those
readers who are not familiar with mathematical standard definitions
and notations such as matrix, group, basis, equivalence class, etc.\ 
are referred to the first paper \cite{KliRRT95} in this series.
However, we tried to make this work as self-contained as possible and
hope that it should be understandable for readers with a rather
limited knowledge of mathematics.

The present version almost fully coincides with \cite{BabCKP97}, 
see Section~\ref{sect:disc} for more details.

\section{Preliminaries}
An {\it undirected graph\/} is a pair $\Gamma=(\Omega,E)$ consisting
of finite sets $\Omega$ and $E$, the {\it vertices\/} and the {\it
  edges\/}.  Each edge connects two different vertices $u$ and $v$
from $\Omega$ and is denoted by $\{u,v\}$. This means that each
element from $E$ is an unordered pair of different vertices from
$\Omega$.  A {\it directed graph\/} is a pair $\Delta=(\Omega,R)$ with
vertex set $\Omega$ and {\it arc set\/} $R$, where each arc, denoted
by $(u,v)$, links two different vertices $u$ and $v$ and additionally
is assigned a direction, namely from $u$ to $v$.  Each element of $R$
is an ordered pair of different vertices from $\Omega$.  If a vertex
$u$ belongs to an edge or an arc then $u$ is said to be {\it
  incident\/} to the edge or arc.  Often it is convenient or useful to
consider an undirected graph $\Gamma$ as a directed graph $\Delta$
with each edge $\{u,v\}$ replaced by two arcs $(u,v)$ and $(v,u)$.

Usually, a directed or undirected graph is given either by its diagram
or by its adjacency matrix.  A {\it diagram\/} is a drawing on the
plane consisting of small circles which represent the vertices and
lines between pairs of vertices which represent the edges.  An arc
$(u,v)$ is indicated by an arrow starting in vertex $u$ and ending in
vertex $v$.  A more abstract representation is the {\it
  (0,1)-adjacency matrix\/} $A=(a_{uv})$.  In order to make evident
that the matrix $A$ represents a graph $\Gamma$, we will also write
$A(\Gamma)$.  The rows and columns of $A$ are indexed by the elements
of $\Omega$, which for sake of simplicity are often numbered by
$1,2,\ldots,n$ with $n=|\Omega|$ (or, as for example in the
computer package \COCO, see below, by numbers $0,1,\ldots, n-1$). Thus,
$A$ is a $n \times n-$matrix. The entry $a_{uv}$ is equal to $1$ if
the edge $\{u,v\}$, respectively the arc $(u,v)$, exists and $0$
otherwise.  Note that the (0,1)-adjacency matrix of an undirected
graph is symmetric with respect to the main diagonal, whereas in
general this is not the case for directed graphs.

Sometimes it is necessary to deal with (undirected or directed)
multigraphs.  In a {\it multigraph\/} each pair of vertices may be
connected by more than one edge or arc. The number of edges resp.\ 
arcs between two vertices is called the {\it multiplicity\/} of the
edge resp.\ arc. In the diagram multiple edges are drawn as parallel
lines, multiple arcs as parallel arrows, in the {\it adjacency
  matrix\/} the entry $a_{uv}$ denotes the multiplicity of the edge or
arc connecting $u$ and $v$.

The {\it complete directed graph\/} is the graph with $n$ vertices
where all $n(n-1)$ arcs are present.

For certain purposes it is more convenient to consider graphs with
{\it loops\/}.  A loop is an arc connecting a vertex with itself.  In
this sense, a complete directed graph with loops has $n^2$ arcs.  In
particular, there is one additional arc $(u,u)$ for each vertex $u$.
The main advantage is that vertices can be identified with the
corresponding loops, which considerably simplifies our notation.

The most general notion of a graph is the colored graph.

In a {\it colored\/} complete directed graph $\Delta$ all vertices and
all arcs are assigned colors in such a way that the colors of the
vertices are different from the colors of the arcs. Assume that
$\{0,1,\ldots,s-1\}$ are the colors of the vertices and let $\Omega_j$
denote the vertices which are assigned color $j$.  Then
$\Omega=\Omega_0 \cup \Omega_1 \cup \ldots \cup \Omega_{s-1}$ is a
partition of the vertex set of $\Delta$.  Similarly, if
$\{s,s+1,\ldots,r-1\}$ are the colors of the arcs and $R_k$ denotes
the arcs of color $k$ then $R=R_s \cup R_{s+1} \cup \ldots \cup R_{r-1}$
is a partition of the arc set of $\Delta$.  Each colored complete
directed graph can be represented by its adjacency matrix $A=(a_{uv})$
which contains in the $u$th row and $v$th column the color of the arc
$(u,v)$, that means $a_{uv}=k$ if and only if $(u,v) \in R_k$. The
entry in the $u$th row and $u$th column is the color of the vertex
$u$, thus $a_{uu}=j$ if and only if $u \in \Omega_j$.

Obviously, any undirected or directed graph can be considered as a colored complete 
directed graph with three colors. The vertices are assigned color $0$, the edges 
(arcs) and nonedges (nonarcs) are assigned colors $1$ and $2$. In the case of a 
multigraph, different colors of arcs correspond to different multiplicities.
In this sense, any chemical structure can be seen as a colored complete graph. 
The colors $a_{uu}$ can be interpreted as modes of vertices, for example names of
atoms in a molecular graph, the colors $a_{uv}$ reflect the multiplicity of bonds
or symbolize certain chains of atoms. 
Figure 1 shows a chemical structure and the adjacency matrix $A$ of the corresponding 
colored complete directed graph (interesting properties of this
compound are discussed in \cite{DunB95})
is given below;
here 0, 1, 2, 3 stands for the
atoms of C, N, O, H respectively, 4 denotes usual bond, 5 double bond,
all other pairs of atoms are denoted by 6, upper superscripts denote
the labels of atoms.

\begin{center}
\setlength{\unitlength}{0.00083300in}%
\begingroup\makeatletter\ifx\SetFigFont\undefined%
\gdef\SetFigFont#1#2#3#4#5{%
  \reset@font\fontsize{#1}{#2pt}%
  \fontfamily{#3}\fontseries{#4}\fontshape{#5}%
  \selectfont}%
\fi\endgroup%
\begin{picture}(4200,5295)(6901,-6100)
\thicklines
\put(8216,-4082){\line( 1,-1){675}}
\put(8263,-4037){\line( 1,-1){675}}
\put(8340,-2828){\line( 1, 1){592}}
\put(8381,-2868){\line( 1, 1){592.500}}
\put(9151,-2236){\line( 1,-1){675}}
\put(9826,-4036){\line(-1,-1){600}}
\put(8101,-3736){\line( 0, 1){600}}
\put(9976,-3136){\line( 0,-1){600}}
\put(9901,-3136){\line( 0,-1){600}}
\put(10051,-4036){\line( 1,-1){375}}
\put(10787,-4561){\line( 1, 0){239}}
\put(10779,-4486){\line( 1, 0){247}}
\put(10501,-4636){\line( 0,-1){300}}
\put(10576,-5011){\makebox(6.6667,10.0000){\SetFigFont{10}{12}{\rmdefault}{\mddefault}{\updefault}.}}
\put(10576,-4636){\line( 0,-1){300}}
\put(9049,-4970){\line( 0,-1){255}}
\put(9049,-5563){\line( 0,-1){255}}
\put(8026,-4036){\line(-1,-1){324}}
\put(7426,-4561){\line(-1, 0){331}}
\put(7426,-4486){\line(-1, 0){324}}
\put(7501,-4636){\line( 0,-1){317}}
\put(7576,-4636){\line( 0,-1){317}}
\put(10159,-2803){\line( 1, 1){267}}
\put(9076,-1945){\line( 0, 1){309}}
\put(8926,-1411){\line(-1, 1){375}}
\put(8977,-1351){\line(-1, 1){375}}
\put(8026,-2911){\line(-1, 1){375}}
\put(9226,-1336){\line( 1, 1){300}}
\put(9274,-1392){\line( 1, 1){300}}
\put(9001,-2161){\makebox(0,0)[lb]{\smash{\SetFigFont{12}{14.4}{\rmdefault}{\mddefault}{\updefault}$C^1$}}}
\put(9901,-3061){\makebox(0,0)[lb]{\smash{\SetFigFont{12}{14.4}{\rmdefault}{\mddefault}{\updefault}$C^2$}}}
\put(9901,-3961){\makebox(0,0)[lb]{\smash{\SetFigFont{12}{14.4}{\rmdefault}{\mddefault}{\updefault}$C^3$}}}
\put(9001,-4861){\makebox(0,0)[lb]{\smash{\SetFigFont{12}{14.4}{\rmdefault}{\mddefault}{\updefault}$C^4$}}}
\put(8101,-3961){\makebox(0,0)[lb]{\smash{\SetFigFont{12}{14.4}{\rmdefault}{\mddefault}{\updefault}$C^5$}}}
\put(8101,-3061){\makebox(0,0)[lb]{\smash{\SetFigFont{12}{14.4}{\rmdefault}{\mddefault}{\updefault}$C^6$}}}
\put(9001,-1561){\makebox(0,0)[lb]{\smash{\SetFigFont{12}{14.4}{\rmdefault}{\mddefault}{\updefault}$N^9$}}}
\put(9601,-961){\makebox(0,0)[lb]{\smash{\SetFigFont{12}{14.4}{\rmdefault}{\mddefault}{\updefault}$O^{15}$}}}
\put(8401,-961){\makebox(0,0)[lb]{\smash{\SetFigFont{12}{14.4}{\rmdefault}{\mddefault}{\updefault}$O^{14}$}}}
\put(10501,-2461){\makebox(0,0)[lb]{\smash{\SetFigFont{12}{14.4}{\rmdefault}{\mddefault}{\updefault}$H^{18}$}}}
\put(10501,-4561){\makebox(0,0)[lb]{\smash{\SetFigFont{12}{14.4}{\rmdefault}{\mddefault}{\updefault}$N^7$}}}
\put(11101,-4561){\makebox(0,0)[lb]{\smash{\SetFigFont{12}{14.4}{\rmdefault}{\mddefault}{\updefault}$O^{10}$}}}
\put(10501,-5161){\makebox(0,0)[lb]{\smash{\SetFigFont{12}{14.4}{\rmdefault}{\mddefault}{\updefault}$O^{11}$}}}
\put(9001,-5461){\makebox(0,0)[lb]{\smash{\SetFigFont{12}{14.4}{\rmdefault}{\mddefault}{\updefault}$O^{16}$}}}
\put(9001,-6061){\makebox(0,0)[lb]{\smash{\SetFigFont{12}{14.4}{\rmdefault}{\mddefault}{\updefault}$H^{19}$}}}
\put(7501,-4561){\makebox(0,0)[lb]{\smash{\SetFigFont{12}{14.4}{\rmdefault}{\mddefault}{\updefault}$N^8$}}}
\put(6901,-4561){\makebox(0,0)[lb]{\smash{\SetFigFont{12}{14.4}{\rmdefault}{\mddefault}{\updefault}$O^{12}$}}}
\put(7501,-5161){\makebox(0,0)[lb]{\smash{\SetFigFont{12}{14.4}{\rmdefault}{\mddefault}{\updefault}$O^{13}$}}}
\put(7501,-2461){\makebox(0,0)[lb]{\smash{\SetFigFont{12}{14.4}{\rmdefault}{\mddefault}{\updefault}$H^{17}$}}}
\end{picture}

\vspace*{0.3cm}

Figure 1  
\end{center}

\[
A=\left(
  \begin{array}{*{19}{c}}
    0 & 4 & 6 & 6 & 6 & 5 & 6 & 6 & 4 & 6 & 6 & 6 & 6 & 6 & 6 & 6 & 6 & 6 & 6\\
    4 & 0 & 5 & 6 & 6 & 6 & 6 & 6 & 6 & 6 & 6 & 6 & 6 & 6 & 6 & 6 & 6 & 4 & 6\\
    6 & 5 & 0 & 4 & 6 & 6 & 4 & 6 & 6 & 6 & 6 & 6 & 6 & 6 & 6 & 6 & 6 & 6 & 6\\
    6 & 6 & 4 & 0 & 5 & 6 & 6 & 6 & 6 & 6 & 6 & 6 & 6 & 6 & 6 & 4 & 6 & 6 & 6\\
    6 & 6 & 6 & 5 & 0 & 4 & 6 & 4 & 6 & 6 & 6 & 6 & 6 & 6 & 6 & 6 & 6 & 6 & 6\\
    5 & 6 & 6 & 6 & 4 & 0 & 6 & 6 & 6 & 6 & 6 & 6 & 6 & 6 & 6 & 6 & 4 & 6 & 6\\
    6 & 6 & 4 & 6 & 6 & 6 & 1 & 6 & 6 & 5 & 5 & 6 & 6 & 6 & 6 & 6 & 6 & 6 & 6\\
    6 & 6 & 6 & 6 & 4 & 6 & 6 & 1 & 6 & 6 & 6 & 5 & 5 & 6 & 6 & 6 & 6 & 6 & 6\\
    4 & 6 & 6 & 6 & 6 & 6 & 6 & 6 & 1 & 6 & 6 & 6 & 6 & 5 & 5 & 6 & 6 & 6 & 6\\
    6 & 6 & 6 & 6 & 6 & 6 & 5 & 6 & 6 & 2 & 6 & 6 & 6 & 6 & 6 & 6 & 6 & 6 & 6\\
    6 & 6 & 6 & 6 & 6 & 6 & 5 & 6 & 6 & 6 & 2 & 6 & 6 & 6 & 6 & 6 & 6 & 6 & 6\\
    6 & 6 & 6 & 6 & 6 & 6 & 6 & 5 & 6 & 6 & 6 & 2 & 6 & 6 & 6 & 6 & 6 & 6 & 6\\
    6 & 6 & 6 & 6 & 6 & 6 & 6 & 5 & 6 & 6 & 6 & 6 & 2 & 6 & 6 & 6 & 6 & 6 & 6\\
    6 & 6 & 6 & 6 & 6 & 6 & 6 & 6 & 5 & 6 & 6 & 6 & 6 & 2 & 6 & 6 & 6 & 6 & 6\\
    6 & 6 & 6 & 6 & 6 & 6 & 6 & 6 & 5 & 6 & 6 & 6 & 6 & 6 & 2 & 6 & 6 & 6 & 6\\
    6 & 6 & 6 & 4 & 6 & 6 & 6 & 6 & 6 & 6 & 6 & 6 & 6 & 6 & 6 & 2 & 6 & 6 & 4\\
    6 & 6 & 6 & 6 & 6 & 4 & 6 & 6 & 6 & 6 & 6 & 6 & 6 & 6 & 6 & 6 & 3 & 6 & 6\\
    6 & 4 & 6 & 6 & 6 & 6 & 6 & 6 & 6 & 6 & 6 & 6 & 6 & 6 & 6 & 6 & 6 & 3 & 6\\
    6 & 6 & 6 & 6 & 6 & 6 & 6 & 6 & 6 & 6 & 6 & 6 & 6 & 6 & 6 & 4 & 6 & 6 & 3
  \end{array}
  \right)
\]

A {\it permutation\/} $f$ acting on a finite set $\Omega$ is a bijective mapping from 
$\Omega$ onto itself. For each permutation $f$ we denote by $v=u^f$ the image $v$ of 
an element $u \in \Omega$. Let $S_n$ be the {\it symmetric group\/} of {\it degree\/} 
$n$, i.e.\ the group of all permutations acting on the set $\Omega$ with $n=|\Omega|$. 
Each subgroup $G$ of $S_n$ is called a {\it permutation group\/} of degree $n$. The 
notation $(G,\Omega)$ indicates that the permutation group $G$ acts on the set $\Omega$.

An {\it automorphism\/} of a colored complete directed graph $\Delta=(\Omega,R)$ is a 
permutation $g$ on $\Omega$ which preserves the colors of the vertices and arcs, 
i.e.\ which fulfills $u \in \Omega_j \Leftrightarrow u^g \in \Omega_j$ 
and $(u,v) \in R_k \Leftrightarrow (u^g,v^g) \in R_k$ 
for all $u,v \in \Omega$ and all colors $j,k$. 
It is easy to realize that the set of all automorphisms of a graph $\Delta$ forms 
a group. This group is called the {\it automorphism group\/} of $\Delta$ and is 
denoted by $Aut(\Delta)$. Clearly, $G=Aut(\Delta)$ is a permutation group acting 
on $\Omega$.

Let $(G,\Omega)$ be a permutation group acting on $\Omega$. We define a binary relation 
$\approx$ on $\Omega$ in such a way that $u \approx v$ holds for two elements 
$u,v \in \Omega$ if and only if there exists a permutation $g \in G$ with $v=u^g$. 
Since $(G,\Omega)$ is a group, the relation $\approx$ is an equivalence relation on 
$\Omega$, its equivalence classes are called {\it orbits\/} (or {\it 1-orbits\/}) of 
$(G,\Omega)$. The set {\it 1-orb(G,$\Omega$)\/} of the orbits of $(G,\Omega)$ forms a 
partition of the set $\Omega$.

The {\it automorphism partition\/} of the vertex set $\Omega$ of a colored complete 
directed graph $\Delta$ is the set {\it 1-orb(Aut($\Delta$),$\Omega$)\/} of the orbits 
of its automorphism group. Obviously, if two vertices $u$ and $v$ belong to the same 
1-orbit, then there is an automorphism $g$ which maps $u$ onto $v$.
The {\it automorphism partitioning problem\/} is the problem of finding the automorphism 
partition of a graph.

To give an example, the bijective mapping $g$ on
$\Omega=\{1,2,\ldots,19\}$ defined by
$$(1)(2,6)(3,5)(4)(7,8)(9)(10,12)(11,13)(14,15)(16)(17,18)(19)$$ is an
automorphism of the colored complete directed graph $\Delta$ which
belongs to the structure in Figure 1. The automorphism partition of
$\Delta$ is $$\left\{\{1\}, \{2,6\}, \{3,5\}, \{4\}, \{7,8\}, \{9\},
  \{10,11,12,13\}, \{14,15\}, \{16\}, \{17,18\}, \{19\}\right\}.$$

A graph $\Delta$ is called {\it vertex-transitive\/} if for any two vertices $u$ 
and $v$ there exists at least one automorphism such that $v=u^g$. Obviously, if 
a graph is vertex-transitive then its automorphism partition is trivial, meaning 
that there is exactly one orbit containing all vertices from $\Omega$.

Commonly, a permutation $f$ on $\Omega$ is represented by a so called
{\it permutation matrix\/} $M(f)=(m_{uv})$. This $n \times n-$matrix has entries $0$ 
and $1$ with $m_{uv}=1$ if and only if $v=u^f$. It is easy to see that a permutation 
matrix has exactly one entry equal to $1$ in every row and in every column, all other 
entries are $0$. In fact, this property is a necessary and sufficient condition for a 
matrix to be a permutation matrix.
Now the property of a permutation to be an automorphism of a graph can be reformulated
in terms of matrices. Namely, a permutation matrix $M$ determines an automorphism 
$g$ of $\Delta$ if and only if $M$ commutes with the adjacency matrix $A$ of 
$\Delta$. This means that the equality
\[ M \cdot A=A \cdot M \]
holds. For example, let $\Gamma$ be the undirected graph depicted in Figure 2.

\unitlength0.5cm
\hspace*{4.9cm}
\begin{picture}(15,6.5)
\thicklines
\put(2.1,2.5){\circle*{0.5}}
\put(9.9,2.5){\circle*{0.5}}
\put(6,0.5){\circle*{0.5}}
\put(6,4.5){\circle*{0.5}}

\put(6,0.75){\line(0,1){4}}
\put(2.3,2.7){\line(2,1){3.5}}
\put(2.3,2.3){\line(2,-1){3.5}}
\put(9.7,2.7){\line(-2,1){3.5}}
\put(9.7,2.3){\line(-2,-1){3.5}}

\put(1,2.25){$1$}
\put(5.8,5.2){$2$}
\put(10.7,2.25){$3$}
\put(5.8,-0.85){$4$}

\end{picture}

\vspace*{0.3cm}

\begin{center} Figure 2 \end{center}

Then
\[A=\left( 
\begin{array}{cccc}
0 & 1 & 2 & 1 \\
1 & 0 & 1 & 1 \\
2 & 1 & 0 & 1 \\
1 & 1 & 1 & 0              
\end{array} \right) \]
is the adjacency matrix of the corresponding colored complete directed graph $\Delta$.
The permutation $g$ of $\Omega=\{1,2,3,4\}$ defined by $1^g=3$, $2^g=4$, $3^g=1$ and 
$4^g=2$ is an automorphism of $\Delta$ with associated permutation matrix
\[M=  \left( 
\begin{array}{cccc}
0 & 0 & 1 & 0 \\
0 & 0 & 0 & 1 \\
1 & 0 & 0 & 0 \\
0 & 1 & 0 & 0              
\end{array} \right).\]

Given a permutation group $(G,\Omega)$, let us now consider the set of permutation 
matrices $M(G)=\{M(g)\,|\,g \in G\}$. A graph $\Delta$ is called {\it invariant\/} 
with respect to the permutation group $(G,\Omega)$ if and only if its adjacency matrix 
commutes with all permutation matrices from $M(G)$.
Let us further consider the set $V(G,\Omega)$ of all $n \times n$-matrices $B$ which 
commute with matrices from $M(G)$, i.e.
\[V(G,\Omega)=\{B\,\,|\,\,M(g) \cdot B = B \cdot M(g) \quad \mbox{for all} 
\quad g \in G\}.\]
$V(G,\Omega)$ is called the {\it centralizer algebra\/} of the permutation group 
$(G,\Omega)$ (the notation $V(G,\Omega)$ stems from the German word
``Vertauschungsring'', the use of which goes back to I. Schur and H. Wielandt).
It is easy to realize that the set of nonnegative integer matrices $B$
from $V(G,\Omega)$ coincides with the set of adjacency matrices of multigraphs which 
are invariant with respect to $(G,\Omega)$. A centralizer algebra $V(G,\Omega)$ is 
known to have the following main properties:

(i) \hspace*{0.62cm} $V(G,\Omega)$ can be considered as a linear space with basis 
$A_0,A_1,\ldots,A_{r-1}$, \\
\hspace*{1.15cm} where each $A_i$ is a $(0,1)$-matrix.  \vspace*{0.2cm} 

(ii) \hspace*{0.5cm} $\sum_{i=0}^{r-1} A_i = J$, where $J$ is the matrix with all 
entries equal to $1$.  \vspace*{0.2cm} 

(iii) \hspace*{0.4cm} For every matrix $A_i$ there exists a matrix $A_j$ with
$A_i^t=A_j$,\\
\hspace*{1.15cm} where $A_i^t$ denotes the transposed matrix of $A_i$.
\vspace*{0.2cm} 

To each permutation group $(G,\Omega)$ we can associate a new induced permutation group
$(G,\Omega^2)$, where for $f \in G$ and $(u,v) \in \Omega^2$ we define\\
\hspace*{7cm} $(u,v)^f = (u^f,v^f)$.\\
Let {\it 2-orb(G,$\Omega$)\/} be the set of orbits of the induced action of $G$ on 
$\Omega^2$.
This partition of the set of all ordered pairs of elements of $\Omega$ is called the 
partition into {\it 2-orbits\/} of $(G,\Omega)$. Each member $R_i\in
\twoo(G,\Omega)$ of this partition defines a graph
$\Gamma_i=\left(\Omega, R_i\right)$ with the adjacency matrix
$A_i=A\left(\Gamma_i\right)$. It turns out that the matrices $A_0,
A_1, \ldots, A_{r-1}$ mentioned in (i) coincide with the latter
adjacency matrices.

More precisely, the basis matrices $A_0,A_1,\ldots,A_{r-1}$ of the centralizer algebra 
$V(G,\Omega)$ correspond to the $2-$orbits of the permutation group $(G,\Omega)$ 
in the following manner.
The $(u,v)-$entry of the basis matrix $A_i$ is equal to $1$ if and only if $(u,v)$ 
belongs to the $i$-$th$ $2-$orbit. In this sense, the set of all $2-$orbits can be 
represented very conveniently by a single matrix of the form
$A = \sum_{i=0}^{r-1} i\cdot A_i$. This means that two pairs $(u,v)$ and $(u',v')$
belong to the same $2-$orbit if and only if the corresponding entries in the matrix 
$A$ are equal.

Let us consider as an example the set $$G=\left\{(1)(2)(3)(4)(5),
  (1,2,3), (1,3,2), (1,2)(4,5), (1,3)(4,5), (2,3)(4,5)\right\}$$ of
  permutations acting on the set $\Omega=\{1,2,3,4,5\}$. It can be
  easily proved that $(G,\Omega)$ is a permutation group (cf. 4.16 in
  \cite{KliRRT95}). Now consider the corresponding set of
  permutation matrices
\[M(G)=\{
\left( 
\begin{array}{ccccc}
1 & 0 & 0 & 0 & 0 \\
0 & 1 & 0 & 0 & 0 \\
0 & 0 & 1 & 0 & 0 \\
0 & 0 & 0 & 1 & 0 \\
0 & 0 & 0 & 0 & 1             
\end{array} \right),
\left(
\begin{array}{ccccc}
0 & 1 & 0 & 0 & 0 \\
0 & 0 & 1 & 0 & 0 \\
1 & 0 & 0 & 0 & 0 \\
0 & 0 & 0 & 1 & 0 \\
0 & 0 & 0 & 0 & 1             
\end{array} \right),
\left(
\begin{array}{ccccc}
0 & 0 & 1 & 0 & 0 \\
1 & 0 & 0 & 0 & 0 \\
0 & 1 & 0 & 0 & 0 \\
0 & 0 & 0 & 1 & 0 \\
0 & 0 & 0 & 0 & 1             
\end{array} \right),\]
\[
\left(
\begin{array}{ccccc}
0 & 1 & 0 & 0 & 0 \\
1 & 0 & 0 & 0 & 0 \\
0 & 0 & 1 & 0 & 0 \\
0 & 0 & 0 & 0 & 1 \\
0 & 0 & 0 & 1 & 0             
\end{array} \right),
\left(
\begin{array}{ccccc}
0 & 0 & 1 & 0 & 0 \\
0 & 1 & 0 & 0 & 0 \\
1 & 0 & 0 & 0 & 0 \\
0 & 0 & 0 & 0 & 1 \\
0 & 0 & 0 & 1 & 0             
\end{array} \right),
\left(
\begin{array}{ccccc}
1 & 0 & 0 & 0 & 0 \\
0 & 0 & 1 & 0 & 0 \\
0 & 1 & 0 & 0 & 0 \\
0 & 0 & 0 & 0 & 1 \\
0 & 0 & 0 & 1 & 0             
\end{array} \right)
\}.\]

It turns out (cf. 5.4 in \cite{KliRRT95}) that the centralizer algebra
$V(G,\Omega)$ has the basis
\[A_0=\left( 
\begin{array}{ccccc}
1 & 0 & 0 & 0 & 0 \\
0 & 1 & 0 & 0 & 0 \\
0 & 0 & 1 & 0 & 0 \\
0 & 0 & 0 & 0 & 0 \\
0 & 0 & 0 & 0 & 0             
\end{array} \right),
A_1=\left( 
\begin{array}{ccccc}
0 & 0 & 0 & 0 & 0 \\
0 & 0 & 0 & 0 & 0 \\
0 & 0 & 0 & 0 & 0 \\
0 & 0 & 0 & 1 & 0 \\
0 & 0 & 0 & 0 & 1             
\end{array} \right),
A_2=\left( 
\begin{array}{ccccc}
0 & 1 & 1 & 0 & 0 \\
1 & 0 & 1 & 0 & 0 \\
1 & 1 & 0 & 0 & 0 \\
0 & 0 & 0 & 0 & 0 \\
0 & 0 & 0 & 0 & 0             
\end{array} \right),\]
\[A_3=\left( 
\begin{array}{ccccc}
0 & 0 & 0 & 0 & 0 \\
0 & 0 & 0 & 0 & 0 \\
0 & 0 & 0 & 0 & 0 \\
0 & 0 & 0 & 0 & 1 \\
0 & 0 & 0 & 1 & 0             
\end{array} \right),
A_4=\left( 
\begin{array}{ccccc}
0 & 0 & 0 & 1 & 1 \\
0 & 0 & 0 & 1 & 1 \\
0 & 0 & 0 & 1 & 1 \\
0 & 0 & 0 & 0 & 0 \\
0 & 0 & 0 & 0 & 0          
\end{array} \right),
A_5=\left( 
\begin{array}{ccccc}
0 & 0 & 0 & 0 & 0 \\
0 & 0 & 0 & 0 & 0 \\
0 & 0 & 0 & 0 & 0 \\
1 & 1 & 1 & 0 & 0 \\
1 & 1 & 1 & 0 & 0             
\end{array} \right),\]
thus we obtain
\[A=\left( 
\begin{array}{ccccc}
0 & 2 & 2 & 4 & 4 \\
2 & 0 & 2 & 4 & 4 \\
2 & 2 & 0 & 4 & 4 \\
5 & 5 & 5 & 1 & 3 \\
5 & 5 & 5 & 3 & 1             
\end{array} \right).\]

To our knowledge, B.Yu. Weisfeiler and A.A. Leman were the first to consider
the more general problem (in comparison with the automorphism partition) of finding
the set of $2-$orbits of $(Aut(\Delta),\Omega)$ for a given graph $\Delta$. In this 
setting the automorphism partition of a graph is a byproduct of the determination of 
the $2-$orbits: the $1-$orbit of a vertex $u$ is simply the $2-$orbit of the pair
$(u,u)$
(for brevity, we will sometimes speak about the orbits of a graph, of a vertex, etc., 
instead of the orbits of the automorphism group of the graph).

In the following it will be our goal to determine the $2-$orbits of a given graph or,
equivalently, to find the $2-$orbit matrix $A$.

\section{Stabilization Procedures}
First attempts to attack the automorphism partitioning problem date
back approximately thirty years.  All these approaches try to find the
$1-$orbits of a given graph.  Usually, the classical paper
\cite{Mor65} by H.L.\ Morgan is considered to be the first procedure
for graph stabilization. Given an undirected graph
$\Gamma=(\Omega,E)$, the idea is to start with a partition of the
vertex set $\Omega$ according to the valencies.  The {\it valency\/}
of a vertex $u$ is the number of edges which are incident to $u$.  Two
vertices $u$ and $v$ are put into the same class of the partition if
and only if they have equal valencies.  Then this partition is refined
iteratively using the extended valencies.  The {\it extended
  valency\/} of $u$ is defined as the sum of the previous extended
valencies of all {\it neighbours\/} of $u$, i.e.\ of all vertices
which are connected with $u$ by an edge.  Again two vertices are put
into the same class of the partition if and only if they have equal
extended valencies. This iteration procedure terminates as soon as the
{\it stable partition\/} is reached, that is the next partition
coincides with the previous one.

Years later it has been recognized that Morgan's approach is just a special case of 
{\it stabilization of depth 2\/}. This technique works as follows. Assume that we
have a partition $\Omega_0,\Omega_1,\ldots,\Omega_{s-1}$ of the vertex set $\Omega$ 
of $\Gamma$ according to the valencies. Assume further that this partition is
numbered such that $\Omega_0$ contains the vertices of smallest valency, $\Omega_1$
the vertices of second smallest valency, etc. Now for each vertex $u \in \Omega$
we compute a list of length $s$ which contains the valencies of $u$ with respect to each 
class $\Omega_j$ (that means the number of edges connecting $u$ with vertices from 
$\Omega_j$), $j=0,1,\ldots,s-1$. Each class $\Omega_j$ may now be partitioned into 
subclasses, each consisting of vertices with equal lists. In this way we may eventually 
obtain a refinement of the original partition. If this is the case, then the subclasses
are numbered according to the lexicographical ordering of the corresponding lists.
Then we restart the same proceeding with the refined partition. 
If in each class the lists of the vertices are identical, then no further refinement 
is obtained and the procedure stops.

The resulting partition of $\Omega$ is commonly called the {\it total degree 
partition\/} of the graph $\Gamma$ (see e.g.\ \cite{Tin86}).
It is the coarsest partition 
$\Omega_0,\Omega_1,\ldots,\Omega_{s-1}$ of $\Omega$ with the property that every 
two vertices belonging to the same cell $\Omega_j$ have the 
same valencies with respect to any other cell $\Omega_k$, $k=j$
included (the coarsest equitable partition of $\Gamma$ in the sense of
\cite{God93}).

The total degree partition cannot be finer than the automorphism partition, since, 
obviously, a necessary condition for two vertices $u$ and $v$ to belong to the same 
$1-$orbit of $Aut(\Gamma)$ is that they belong to the same class of the total degree 
partition. In fact, every class of the total degree partition is a union of $1-$orbits.
Figure 3 shows a graph where the total degree partition consists of one class only, 
namely $\Omega$, but which is not vertex-transitive, i.e.\ the automorphism partition 
consists of more than one $1-$orbit. This example makes evident the weakness
of the above kind of stabilization. The method will not give any partition for 
{\it regular\/} graphs (graphs where all vertices have the same valencies), not even 
in the case when $Aut(\Gamma)$ is the trivial group consisting of the identity only
(what means that each $1-$orbit consists of a single vertex).

The reader will recognize in Figure 3 {\it ``cuneane''}, cf. 4.26 and
5.8 in \cite{KliRRT95}.

\begin{center}
\setlength{\unitlength}{0.00083300in}%
\begingroup\makeatletter\ifx\SetFigFont\undefined%
\gdef\SetFigFont#1#2#3#4#5{%
  \reset@font\fontsize{#1}{#2pt}%
  \fontfamily{#3}\fontseries{#4}\fontshape{#5}%
  \selectfont}%
\fi\endgroup%
\begin{picture}(3176,2680)(2608,-7884)
\thicklines
\put(4201,-6061){\circle*{150}}
\put(4801,-5461){\circle*{150}}
\put(5701,-5761){\circle*{150}}
\put(5101,-6361){\circle*{150}}
\put(4801,-7561){\circle*{150}}
\put(3301,-7561){\circle*{150}}
\put(2701,-6361){\circle*{150}}
\put(3301,-5761){\circle*{150}}
\put(3301,-5761){\line( 5, 1){1500}}
\put(4801,-5461){\line( 3,-1){900}}
\put(5701,-5761){\line(-1,-1){600}}
\put(5101,-6361){\line(-3, 1){900}}
\put(4201,-6061){\line(-5,-1){1500}}
\put(2701,-6361){\line( 1, 1){600}}
\put(3301,-5761){\line( 0,-1){1800}}
\put(3301,-7561){\line( 1, 0){1500}}
\put(4801,-7561){\line( 1, 4){300}}
\put(4801,-7561){\line( 1, 2){900}}
\put(2701,-6361){\line( 1,-2){600}}
\put(4201,-6061){\line( 1, 1){600}}
\put(3223,-5607){\makebox(0,0)[lb]{\smash{\SetFigFont{12}{14.4}{\rmdefault}{\mddefault}{\updefault}7}}}
\put(4748,-5312){\makebox(0,0)[lb]{\smash{\SetFigFont{12}{14.4}{\rmdefault}{\mddefault}{\updefault}1}}}
\put(5653,-5592){\makebox(0,0)[lb]{\smash{\SetFigFont{12}{14.4}{\rmdefault}{\mddefault}{\updefault}2}}}
\put(5013,-6222){\makebox(0,0)[lb]{\smash{\SetFigFont{12}{14.4}{\rmdefault}{\mddefault}{\updefault}3}}}
\put(4103,-5932){\makebox(0,0)[lb]{\smash{\SetFigFont{12}{14.4}{\rmdefault}{\mddefault}{\updefault}8}}}
\put(2608,-6232){\makebox(0,0)[lb]{\smash{\SetFigFont{12}{14.4}{\rmdefault}{\mddefault}{\updefault}6}}}
\put(3233,-7837){\makebox(0,0)[lb]{\smash{\SetFigFont{12}{14.4}{\rmdefault}{\mddefault}{\updefault}5}}}
\put(4778,-7857){\makebox(0,0)[lb]{\smash{\SetFigFont{12}{14.4}{\rmdefault}{\mddefault}{\updefault}4}}}
\end{picture}
  
\vspace{0.3cm}
 
Figure 3
\end{center}

As already mentioned in the introduction, the problem of the
recognition of graph symmetry 
can be made more precise in mathematical language as an automorphism partitioning problem;
the problem of graph identification corresponds to an isomorphism problem.
Let us give a precise mathematical statement of the latter problem.

An {\it isomorphism\/} from a graph $\Gamma=(\Omega,E)$ to a graph $\Gamma'=(\Omega',E')$ 
is a bijective mapping $h$ from $\Omega$ to $\Omega'$ such that $(u,v) \in E$ if and 
only if $(u^h,v^h) \in E'$. If such a mapping exists then $\Gamma$ and $\Gamma'$ are 
called {\it isomorphic\/}. Obviously, if $\Gamma = \Gamma'$ then an isomorphism 
coincides with an automorphism. The {\it isomorphism problem\/} is the problem of 
deciding whether two graphs are isomorphic or not. Similarly as for the automorphism 
problem there is also a matrix formulation for the isomorphism problem. Namely, two 
graphs $\Gamma$ and $\Gamma'$ with adjacency matrices $A$ resp.\ $A'$ are 
isomorphic if and only if there exists a permutation matrix $M$ with 
\[ M \cdot A = A' \cdot M. \]
If this equality is multiplied from the left with the inverse $M^{-1}$ of $M$
then, using the fact that $M^{-1}=M^t$ holds for each permutation matrix $M$, 
we obtain
\[ A = M^t \cdot A' \cdot M.\]
This modified equality can be interpreted in the
following way. The matrix $A$ is obtainable from the matrix $A'$ by 
permuting simultaneously rows and columns. This corresponds just to a renumbering 
of the vertices.

G.\ Tinhofer has found a very interesting algebraic characterization of total degree 
partitions. He first relaxed the notion of an isomorphism between two graphs
$\Gamma$ and $\Gamma'$ using doubly stochastic matrices. A matrix $X$ is 
{\it doubly stochastic\/} if the entries of $X$ are nonnegative and the sum of 
the entries in each row and in each column is equal to $1$. Note that every 
permutation matrix is doubly stochastic. 
Now let again $A$ and $A'$ be the adjacency matrices of $\Gamma$ resp.\ $\Gamma'$.
Then the two graphs are called {\it doubly stochastic isomorphic\/} if and only if
there exists a doubly stochastic matrix $X$ fulfilling the equality 
\[ X \cdot A = A' \cdot X.\]
Of course, two isomorphic graphs are also doubly stochastic isomorphic, however,
the converse direction is not true in general. For example, the graphs $\Gamma$ 
and $\Gamma'$ of Figure 4 are doubly stochastic isomorphic (choose $X=1/6 \cdot J$),
but they are not isomorphic.

\unitlength0.5cm
\hspace*{1cm}
\begin{picture}(15,6.5)
\thicklines
\put(-0.5,2.25){$\Gamma$}
\put(3,2.5){\circle*{0.5}}
\put(5,0.5){\circle*{0.5}}
\put(5,4.5){\circle*{0.5}}
\put(7.5,0.5){\circle*{0.5}}
\put(7.5,4.5){\circle*{0.5}}
\put(9.5,2.5){\circle*{0.5}}

\put(3.2,2.7){\line(1,1){1.8}}
\put(3.2,2.3){\line(1,-1){1.8}}
\put(5.25,4.5){\line(1,0){2.25}}
\put(5.25,0.5){\line(1,0){2.25}}
\put(9.3,2.7){\line(-1,1){1.8}}
\put(9.3,2.3){\line(-1,-1){1.8}}

\put(1.95,2.25){$1$}
\put(4.8,5.2){$2$}
\put(7.3,5.2){$3$}
\put(10.2, 2.25){$4$}
\put(7.3,-0.65){$5$}
\put(4.8,-0.65){$6$}

\put(14.5,2.25){$\Gamma'$}
\put(18,2.5){\circle*{0.5}}
\put(20,0.5){\circle*{0.5}}
\put(20,4.5){\circle*{0.5}}
\put(22.5,0.5){\circle*{0.5}}
\put(22.5,4.5){\circle*{0.5}}
\put(24.5,2.5){\circle*{0.5}}

\put(18.2,2.7){\line(1,1){1.8}}
\put(18.2,2.3){\line(1,-1){1.8}}
\put(20,0.75){\line(0,1){3.5}}
\put(22.5,0.75){\line(0,1){3.5}}
\put(24.3,2.7){\line(-1,1){1.8}}
\put(24.3,2.3){\line(-1,-1){1.8}}

\put(16.95,2.25){$1$}
\put(19.8,5.2){$2$}
\put(22.3,5.2){$3$}
\put(25.2,2.25){$4$}
\put(22.3,-0.65){$5$}
\put(19.8,-0.65){$6$}

\end{picture}

\vspace*{0.3cm}

\begin{center} Figure 4 \end{center}

Tinhofer proved in \cite{Tin86} that two graphs are doubly stochastic isomorphic 
if and only if they have identical total degree partitions.
To be more precise, the total degree partitions $\Omega_0,\Omega_1,\ldots,\Omega_{s-1}$
of $\Gamma=(\Omega,E)$ and $\Omega_0',\Omega_1',\ldots,\Omega_{s'-1}'$
of $\Gamma'=(\Omega',E')$ are identical, if $s=s'$, $|\Omega_j|=|\Omega_j'|$ and,
for each pair of vertices $u \in \Omega_j$ and $u' \in \Omega_j'$, the valency of
$u$ with respect to $\Omega_k$ is equal to the valency of $u'$ with respect to
$\Omega_k'$, $j,k \in \{0,1,\ldots,s-1\}$.

The shortcoming of the total degree partition, as pointed out above, motivated the
construction of more powerful algorithms for graph stabilization. A well known
approach is to apply the refinement technique, which has already been used for the
computation of the total degree partition, by replacing the valencies of the vertices 
by other invariants of the graph. For example, for each vertex $u$ we may count the
number of cycles of a given length which contain $u$, the number and sizes of
cliques containing $u$, etc. However, finding cycles or cliques in a graph is  
an extremely difficult task which, in general, requires time exponential in the size 
of the graph. Therefore, we should use a criterion which is easy to check.
A reasonable approach is to consider not only configurations consisting of two vertices,
i.e.\ the edges and nonedges (which define the valencies), but to examine also
configurations consisting of three vertices, i.e.\ all triples of vertices. This 
procedure is commonly called {\it stabilization of depth $3$\/}. The algorithm of 
Weisfeiler and Leman, which will be formulated and discussed in great detail in 
Sections 5 and 6, is based on that principle.

We have seen above that stabilization of depth 2 can be associated with a certain 
combinatorial object, namely with the total degree partition of the graph under 
consideration. This immediately implies the question whether there is a similar  
object which belongs to stabilization of depth 3. It turns out that such an object
really exists.
Moreover, this object has not only a combinatorial but also a very interesting 
algebraic description. The details are exposed in the next section.

\section{Cellular Algebras}
A {\it matrix algebra\/} of {\it degree\/} $n$ is a set of $n \times n-$matrices 
which is closed under matrix addition, matrix multiplication and multiplication
of a matrix by a scalar,
i.e.\ if $X$ and $Y$ belong to the matrix algebra and $z$ 
is any real number, then also $X + Y$, $X \cdot Y$ and $z \cdot X$ belong to the
matrix algebra.

A {\it cellular\/} (or {\it coherent\/}) {\it algebra\/} is a matrix algebra which
additionally is closed under Schur-Hadamard (=componentwise) multiplication of matrices 
and under matrix transposition, and which contains the identity matrix $I$ and the 
matrix $J$ all entries of which are equal to 1.
We will denote the Schur-Hadamard product of two matrices $X$ and $Y$ by $X \circ Y$.
Thus, if $X=(x_{uv})$ and $Y=(y_{uv})$ then $X \circ Y = (x_{uv} \cdot y_{uv})$.
Each cellular algebra $W$ has a basis $A_0, A_1, \ldots, A_{r-1}$
(basis of the vector space $W$) consisting of $(0,1)$-matrices which is
called {\it standard basis\/} of $W$, $r$ is the {\it rank\/} of $W$.
It is not hard to see that a standard basis $A_0,A_1,\ldots,A_{r-1}$, if
suitably numbered, satisfies the following properties:

(i) \hspace*{0.62cm} $\sum_{i=0}^{r-1} A_i = J$ \vspace*{0.18cm} 

(ii) \hspace*{0.5cm} $\sum_{i=0}^{q-1} A_i = I$ \hspace*{0.18cm}
     for some $q$ with $q \leq r$               \vspace*{0.18cm}

(iii) \hspace*{0.4cm}  $A_i \circ A_j = 0 \; \Leftrightarrow \; i \neq j, \;\;
     i,j \in \{0,1,\ldots,r-1\}$ \vspace*{0.18cm}

(iv) \hspace*{0.44cm} for each $i \in \{0,1,\ldots,r-1\}$ there is a $j \in
      \{0,1,\ldots,r-1\}$ such that $A_i^t = A_j$ \vspace*{0.18cm}

(v) \hspace*{0.53cm} for each pair $i,j \in \{0,1,\ldots,r-1\}$ we have \vspace*{0.18cm}

\hspace*{1.4cm}
$A_i A_j = p_{ij}^0 A_0 + p_{ij}^1 A_1 + \ldots + p_{ij}^{r-1} A_{r-1}$.

A cellular algebra $W$ with standard basis $A_0,A_1,\ldots,A_{r-1}$ can be
represented in a very convenient way using the matrix $A(W) = \sum_{i=0}^{r-1}
i \cdot A_i$, called the {\it adjacency matrix of the cellular algebra\/} $W$.
This matrix is unique up to the numbering of the basis matrices.

The nonnegative integers $p_{ij}^k$ are called the {\it structure constants\/}
of $W$. These numbers have a very nice geometric interpretation.  We consider
the colored complete directed graph $\Delta=(\Omega,R)$ which belongs to the
adjacency matrix of $W$. For our purposes, it is very convenient to identify
each vertex $u$ in $\Delta$ with the corresponding loop $(u,u)$, i.e.\ we deal
with the complete directed graph with loops.  Then the matrices $A_k$
correspond in a natural way to the arc sets $R_k$ of color $k$ (the first $q$
matrices of the standard basis virtually represent the vertices of $\Delta$;
the induced partition of the vertex set is called the {\it standard
partition\/}\footnote[1]{ More rigorously, we should speak about {\it standard
partition of depth 3\/} in order to emphasize that the partition is induced by
stabilization of depth 3.
In the following, we briefly call it the standard partition.} of the graph).
An arc $(u,v)$ has color $k$ if and only if the matrix $A_k$ has entry $1$ in
the $u$th row and $v$th column. The entry $(u,v)$ in the matrix $A_i \cdot A_j$
is the number of directed paths of length $2$ from vertex $u$ to vertex $v$,
such that the first step is an arc of color $i$ and the second step is an arc
of color $j$. The decomposition (v) implies that for any arc $(u,v)$ of a fixed
color $k$ the number of paths of length $2$ from $u$ to $v$, such that the
first step is of color $i$ and the second step is of color $j$, is the same and
equal to $p_{ij}^k$ (see Figure 5).  In fact, if $A_0,A_1,\ldots,A_{r-1}$
fulfill (i)-(iv), then this condition is sufficient for these matrices to be
the standard basis of a cellular algebra.

\unitlength0.5cm
\begin{picture}(15,7)
\thicklines
\put(12,0){\circle*{0.5}}
\put(20,0){\circle*{0.5}}
\put(16,4){\circle*{0.5}}

\put(12.25,0){\vector(1,0){7.5}}
\put(12.2,0.15){\vector(1,1){3.6}}
\put(16.2,3.85){\vector(1,-1){3.6}}

\put(10.7,-0.2){$u$}
\put(21,-0.2){$v$}
\put(15.7,5){$w$}
\put(15.8,-1){$k$}
\put(13.3,2.2){$i$}
\put(18.7,2.2){$j$}

\end{picture}

\vspace*{0.5cm}

\begin{center} Figure 5 \end{center}

It is not hard to see (see \cite{KliRRT95})
that each centralizer algebra is also a cellular algebra
(therefore, the matrix $A$ given at the end of Section 2 represents a cellular 
algebra of rank 6 with standard basis $A_0,A_1,\ldots,A_5$).
Let $W$ be a cellular algebra.
If a group with the centralizer algebra $W$ exists then the cellular algebra $W$ is
called {\it Schurian\/}, after I.\ Schur, who was the
first to investigate cellular algebras (actually he was using a different terminology of
group rings, see \cite{KliRRT95} for details).
The importance of Schurian cellular algebras stems from the fact that, as already
indicated in Section 2, its basis matrices correspond to the $2-$orbits. Moreover,
the diagonal matrices of the basis correspond to the $1-$orbits.

Given any $n \times n-$matrix $X$, the {\it cellular algebra $W(X)$ generated by\/} 
$X$ is defined to be the {\it smallest\/} cellular algebra which contains $X$. 
It is important to know that this definition is rigorous. That means the resulting 
algebra is well defined and unique (for a detailed explanation see \cite{KliRRT95}).
As a consequence, we are able to associate with each graph $\Gamma$ a matrix 
algebra, namely the cellular algebra $W(A)$ which is generated by the adjacency 
matrix $A$ of $\Gamma$. We will also write $W(\Gamma)$ in order to indicate that the
cellular algebra corresponds to the graph $\Gamma$.

There are a number of graph classes whose associated cellular algebras
are Schurian (for example graphs with a simple spectrum \cite{Pon94a}).
For those graphs the automorphism partition can be immediately deduced
from the cellular algebra.  Namely, in this case the automorphism
partition coincides with the standard partition.  Unfortunately, this
is not the case in general.  The simplest counterexamples can be found
among strongly regular graphs.  An undirected graph is called {\it
  strongly regular\/} (see \cite{HesH71}) if the standard basis of the
associated cellular algebra consists only of $3$ matrices.  The first such
examples of non-Schurian cellular algebras of rank 3 were found by H.\ 
Wielandt \cite{Wie64}, L.C.\ Chang \cite{Cha59}, S.S.\ Shrikhande
\cite{Shr59} and G.M.\ Adel'son-Velski\v{\i} et\ al.\ \cite{AdeWLF69}.
All these graphs have rather large automorphism groups.  Later on,
collaborators of Weisfeiler found an example with the identity
automorphism group (see \cite{Wei76}).

Nevertheless, although we cannot guarantee that we get the $1-$orbits and $2-$orbits for
each graph, extended practical experience indicates that the results obtained by 
the cellular algebras are sufficient, in particular for practically
all chemical graphs.
This is confirmed by the computational results which are presented in Section 9. 

The representation of a cellular algebra as a colored complete directed graph and 
the interpretation of the structure constants $p_{ij}^k$ shows that in a 
cellular algebra implicitly all configurations of a graph consisting of three 
vertices are considered. In other words, the cellular algebra is the algebraic 
object which is associated to stabilization of depth $3$. 

At the end of this section, let us give a precise statement of the problem which now
has to be solved. Given a graph $\Gamma$, we actually deal with two closely
related problems.

\begin{problem}
\item Compute the basis $A_0,A_1,\ldots,A_{r-1}$ of the cellular
  algebra $W(\Gamma)$ (or in other words the colored complete graph
  $\Delta$ which is associated to the cellular algebra $W(\Gamma)$).
\item Construct the colored complete graph $\Delta$ with the structure
  constants $p_{ij}^k$.
\end{problem}

In the next section we give an algebraic description of the algorithm of
Weisfeiler-Leman which settles Problem 1. In Section 6 we present a very
illustrative graph theoretical interpretation which solves Problem 2.

\section{Algebraic Description of the Algorithm}

B.\ Weisfeiler and A.\ Leman were the first to show that the cellular algebra of a
graph can be computed in polynomial time. The proposed method, firstly described in 
Russian in the paper \cite{WeiL68}, has as input any matrix $A$ (the adjacency 
matrix of a graph $\Gamma$) and as output a basis of the cellular algebra $W(A)$ 
generated by $A$. The initial description of the algorithm was too sophisticated, 
clearer ones are given in the English written book \cite{Wei76}, and
also in \cite{Fri89} 
and \cite{Hig87}. We will first explain the main features of the algorithm by combining 
all these ideas, illustrate it by an example, and then present a formal description.

The construction of the cellular algebra $W(A)$ proceeds iteratively. We start with 
the adjacency matrix $A=(a_{uv})$ of an undirected or directed graph $\Gamma$
(the diagonal entries are set to be different from the nondiagonal entries).
This matrix can be written in the form $A = \sum_{i=0}^{r-1}i \cdot A_i$ where 
$A_i$ are $(0,1)-$matrices with $(u,v)-$entry equal to $1$ if and only if $a_{uv}=i$. 
At the end of each iteration we obtain a new set of $(0,1)-$matrices 
$A'_0,A'_1,\ldots,A'_{r'-1}$ which fulfills the properties (i)-(iv) of a cellular 
algebra as stated in Section 4, but which may fail property (v).
In particular, this means that $A'_0,A'_1,\ldots,A'_{r'-1}$ is the basis of a linear 
subspace $S$ which is closed under Schur-Hadamard multiplication and transposition 
and which contains the identity matrix $I$ and the all 1 matrix $J$.
Note that this also holds at the beginning of the procedure in case $A$ is symmetric 
and the values of the diagonal entries are different from all other entries.

Initially, the linear subspace $S$ with basis matrices $A_0,A_1,\ldots,A_{r-1}$ 
will not fulfill property (v) of a cellular algebra, i.e.\ $S$ will not be
closed with respect to matrix multiplication. 
It is easy to see that $S$ contains all products of matrices from $S$
if and only if it contains all products $A_i \cdot A_j$ of basis matrices.
Therefore, we consider in each iteration the linear subspace which is generated by 
all these products $A_0 A_0, A_0 A_1, \ldots, A_{r-1} A_{r-1}$
and compute a basis of $(0,1)-$matrices for it. This process is repeated until
it is stable, this means the basis of the actual iteration coincides with the basis 
of the previous iteration (up to the numbering of the basis matrices).
In that case the subspace $S$ is closed under matrix
multiplication, property (v) is satisfied and, consequently, $S$ is a cellular algebra.

A straightforward method to construct the basis $A'_0,A'_1,\ldots,A'_{r'-1}$ of a 
linear subspace $S$ which is generated by some set of matrices 
$\{B_0,B_1,\ldots,B_{p-1}\}$ (in our case this is just the set 
$\{A_0 A_0,$ $ A_0 A_1, \ldots, A_{r-1} A_{r-1}\}$)
and which is closed under Schur-Hadamard multiplication is described in the paper 
\cite{Hig87}. Let $B=B_0$. If $B=0$ then set $B_i=B_{i+1}$ for
$i=0,1,\ldots,p-1$ and repeat the procedure. Assume that $B=(b_{uv}) \neq 0$ and
let $t$ be a nonzero entry of $B$. Set $D=(d_{uv})$ to be the $(0,1)-$matrix
such that $d_{uv}=1$ if and only if $b_{uv}=t$. Let $B=D \circ B_1$. If $B \neq 0$
then let $D$ be a $(0,1)-$matrix as constructed above. Now let $B=D \circ B_2$
and repeat the procedure. If $B=0$ then let $B=D \circ B_3$ and repeat the procedure.
The last $(0,1)-$matrix $D$ will be the first element $A'_0$ in the basis of $S$.
Let $B_i=B_i-B_i \circ A'_0$, $i=0,1,\ldots,p-1$ and repeat the procedure.
In this way we obtain a set $A'_0,A'_1,\ldots,A'_{r'-1}$ of $(0,1)-$matrices fulfilling
property (iii).

This procedure can be formulated in a more compact and convenient way
(which in fact was the original way used by Weisfeiler--Leman) by introducing 
indeterminates $t_i$, $i=0,1,\ldots,r-1$, which can be considered to represent
the different entries of an adjacency matrix. In this sense, the matrix 
$D = \sum_{i=0}^{r-1} t_i A_i$ represents the adjacency matrix. 
Since matrix multiplication is a noncommutative operation, it is important in the 
following to assume that the indeterminants are noncommuting with respect to 
multiplication, i.e.\ $t_it_j \neq t_jt_i$.
Now compute the product $B=D \cdot D = \sum_{i=0}^{r-1} \sum_{j=0}^{r-1} t_it_jA_iA_j$.
Each entry of this matrix is a sum of products $t_it_j$.
In order to obtain the basis of $S$ we replace equal entries in $B$ by new
indeterminates $t'_i$. Now it is not hard to verify that
the $(0,1)-$matrices $A'_0,A'_1,\ldots,A'_{r'-1}$ of the resulting 
matrix $A'=\sum_{i=0}^{r'-1} t'_i A'_i$ are exactly the basis matrices of $S$.

The following very simple example will illustrate this
procedure. Consider the graph $\Gamma$ which is depicted in Figure 6.

\begin{center}
\setlength{\unitlength}{0.00083300in}%
\begingroup\makeatletter\ifx\SetFigFont\undefined%
\gdef\SetFigFont#1#2#3#4#5{%
  \reset@font\fontsize{#1}{#2pt}%
  \fontfamily{#3}\fontseries{#4}\fontshape{#5}%
  \selectfont}%
\fi\endgroup%
\begin{picture}(2100,1365)(4501,-3094)
\thicklines
\put(5326,-2461){\line( 1, 0){600}}
\put(5326,-2386){\line( 1, 0){600}}
\put(6226,-2236){\line( 1, 1){300}}
\put(6526,-2911){\line(-1, 1){375}}
\put(5026,-2311){\line(-1, 1){375}}
\put(5026,-2536){\line(-1,-1){300}}
\put(5101,-2461){\makebox(0,0)[lb]{\smash{\SetFigFont{12}{14.4}{\rmdefault}{\mddefault}{\updefault}$C^1$}}}
\put(6601,-1861){\makebox(0,0)[lb]{\smash{\SetFigFont{12}{14.4}{\rmdefault}{\mddefault}{\updefault}$H^5$}}}
\put(6601,-3061){\makebox(0,0)[lb]{\smash{\SetFigFont{12}{14.4}{\rmdefault}{\mddefault}{\updefault}$H^6$}}}
\put(4501,-1861){\makebox(0,0)[lb]{\smash{\SetFigFont{12}{14.4}{\rmdefault}{\mddefault}{\updefault}$H^4$}}}
\put(4501,-3061){\makebox(0,0)[lb]{\smash{\SetFigFont{12}{14.4}{\rmdefault}{\mddefault}{\updefault}$H^3$}}}
\put(6001,-2461){\makebox(0,0)[lb]{\smash{\SetFigFont{12}{14.4}{\rmdefault}{\mddefault}{\updefault}$C^2$}}}
\end{picture}

\vspace*{0.4cm}

Figure 6
\end{center}

Here superscripts denote the numbers from $\Omega=\{1,2,3,4,5,6\}$
associated to atoms which form the molecule of ethylene. Let $A$ be the
adjacency matrix of the colored graph $\Gamma$ associated to the
molecular graph depicted in Figure 6 (here 0 stands for the carbon atom, 1
for hydrogen atom, 2 for usual bond, 3 for double bond and 4 means
that there is no bond between the corresponding atoms),
\[
A=\left(
  \begin{array}{*6{c}}
    0 & 3 & 2 & 2 & 4 & 4\\
    3 & 0 & 4 & 4 & 2 & 2\\
    2 & 4 & 1 & 4 & 4 & 4\\
    2 & 4 & 4 & 1 & 4 & 4\\
    4 & 2 & 4 & 4 & 1 & 4\\
    4 & 2 & 4 & 4 & 4 & 1
  \end{array}
  \right) .
\]

Then we get that
\[
D=\left(
  \begin{array}{*6{c}}
    t_0 & t_3 & t_2 & t_2 & t_4 & t_4\\
    t_3 & t_0 & t_4 & t_4 & t_2 & t_2\\
    t_2 & t_4 & t_1 & t_4 & t_4 & t_4\\
    t_2 & t_4 & t_4 & t_1 & t_4 & t_4\\
    t_4 & t_2 & t_4 & t_4 & t_1 & t_4\\
    t_4 & t_2 & t_4 & t_4 & t_4 & t_1
  \end{array}
  \right) ,
\]

\[
B=\left(
  \begin{array}{*6{c}}
    x_0 & x_2 & x_3 & x_3 & x_4 & x_4\\
    x_2 & x_0 & x_4 & x_4 & x_3 & x_3\\
    x_5 & x_6 & x_1 & x_7 & x_8 & x_8\\
    x_5 & x_6 & x_7 & x_1 & x_8 & x_8\\
    x_6 & x_5 & x_8 & x_8 & x_1 & x_7\\
    x_6 & x_5 & x_8 & x_8 & x_7 & x_1
  \end{array}
\right),
\]
where
\begin{eqnarray*}
  x_0 & = & t_0^2+ 2t_2^2+ t_3^2+ 2t_4^2,\\
  x_1 & = & t_1^2+ t_2^2+ 4t_4^2,\\
  x_2 & = & t_0t_3+ 2t_2t_4+ t_3t_0+ 2t_4t_2,\\
  x_3 & = & t_0t_2+ t_2t_1+ t_2t_4+ t_3t_4+ 2t_4^2,\\
  x_4 & = & t_0t_4+ 2t_2t_4+ t_3t_2+ t_4t_1+ t_4^2,\\
  x_5 & = & t_1t_2+ t_2t_0+ t_4t_2+ t_4t_3+ 2t_4^2,\\
  x_6 & = & t_1t_4+ t_2t_3+ t_4t_0+ 2t_4t_2+ t_4^2,\\
  x_7 & = & t_1t_4+ t_2^2+ t_4t_1+ 3t_4^2,\\
  x_8 & = & t_1t_4+ t_2t_4+ t_4t_1+ t_4t_2+ 2t_4^2.
\end{eqnarray*}

Now we proceed with the matrix $A'$,
\[
A'=\left(
  \begin{array}{*6{c}}
    0 & 2 & 3 & 3 & 4 & 4\\
    2 & 0 & 4 & 4 & 3 & 3\\
    5 & 6 & 1 & 7 & 8 & 8\\
    5 & 6 & 7 & 1 & 8 & 8\\
    6 & 5 & 8 & 8 & 1 & 7\\    
    6 & 5 & 8 & 8 & 7 & 1\\
  \end{array}
  \right).
\]

We suggest that the reader repeats the process with the matrix $A'$
instead of $A$ and checks that the matrix $A''$ resulting from the
second iteration coincides with $A'$. This means that $A'$ is in fact
the desired result of the stabilization, namely, the adjacency matrix
$A(W(\Gamma))$. The reader can easily find that the order of the
automorphism group $G=\aut(\Gamma)$ of our graph $\Gamma$ is equal to
$8$ and the matrix $A'$ represents the set of 2-orbits of $(G,
\Omega)$. Hence in this case $W(\Gamma)$ is really a Schurian cellular algebra. 

Here is a formal description of the algorithm (which will be denoted 
in the following by the initials of the authors).
\vspace*{0.2cm}

{\bf Algorithm WL}

Input: the adjacency matrix $A=A(\Gamma)=(a_{uv})$ of colored graph $\Gamma$.

Output: a standard basis $A_0,A_1,\ldots,A_{r-1}$ of the cellular algebra
$W(\Gamma)$, or more exactly the adjacency matrix $A(W(\Gamma))$.

(0) Let $\{0,1,\ldots,s-1\}$ be the set of different entries of $A$.\vspace*{0.1cm}\\
\hspace*{0.5cm} For $k=0,1,\ldots,s-1$ do  \vspace*{0.1cm}\\
\hspace*{1.7cm}
Define $A_k=(a(k)_{uv})$ to be the matrix with \vspace*{0.1cm}\\
\hspace*{1.7cm}
$a(k)_{uv}=1$ if $a_{uv}=k$ and $a(k)_{uv}=0$ otherwise. \vspace*{0.1cm}\\
\hspace*{0.5cm}
Let $r := s$.

(1) Let $D=\sum_{k=0}^{r-1} t_k A_k$, \vspace*{0.1cm}\\
\hspace*{1.1cm} 
where $t_0,t_1,\ldots,t_{r-1}$ are distinct indeterminates, \vspace*{0.1cm}\\
\hspace*{1.1cm} which are noncommuting with respect to multiplication.

(2) Compute the matrix product $B=(b_{uv}) = D \cdot D$. \vspace*{0.1cm}\\
\hspace*{1.1cm} Each entry $b_{uv}$ of $B$ is a sum of products $t_i t_j$.

(3) Determine the set $\{d_0,d_1,\ldots,d_{s-1}\}$  \vspace*{0.1cm}\\
\hspace*{1.1cm} of different expressions among the entries $b_{uv}$.

(4) If $s > r$ then \vspace*{0.1cm}\\
\hspace*{1.7cm} 
For $k=0,1,\ldots,s-1$ do  \vspace*{0.1cm}\\
\hspace*{2.3cm}
Define $A_k=(a(k)_{uv})$ to be the matrix with \vspace*{0.1cm}\\
\hspace*{2.3cm}
$a(k)_{uv}=1$ if $b_{uv}=d_k$ and $a(k)_{uv}=0$ otherwise. \vspace*{0.1cm}\\
\hspace*{1.7cm}
$r := s$. Goto (1).

(5) STOP.
\vspace*{0.5cm}

As already indicated in the previous section, it follows from the representation 
of a cellular algebra by a colored complete directed graph that algorithm WL 
implicitly considers all configurations consisting of three vertices in the graph.
There are other methods (see e.g.\ the references in \cite{RueR90b}) for perception 
of graph symmetry which also use configurations of three vertices. The main advantage 
of algorithm WL over them is the fact that the cellular algebra contains exhaustive 
information which is obtainable from subgraphs of at most three vertices. Thus, the 
resulting partitions of the vertices and pairs of vertices are the finest among
the partitions which can be deduced using these configurations.

An interesting method which is rather close to algorithm WL is described in
\cite{RueR90b}.
It is based on raising the adjacency matrix to its higher powers, evaluating the
entries and partitioning the atoms and pairs of atoms into equivalence classes.
Indeed, this procedure considers subgraphs consisting of three vertices. However, 
some information can get lost when the resulting partitions are coarser than 
those obtained by algorithm WL. In particular, an ``orientation'' of edges cannot occur.
Note that this is possible even if the algorithm WL is applied to an undirected graph.
Consider for example the graph $\Gamma$ in Figure $7$.

\vspace*{1cm}

\unitlength0.5cm
\hspace*{0.9cm}
\begin{picture}(15,3)
\thicklines
\put(3,0){\circle*{0.5}}
\put(6,0){\circle*{0.5}}
\put(8.5,1.65){\circle*{0.5}}
\put(8.5,-1.65){\circle*{0.5}}

\put(3.25,0){\line(1,0){2.5}}
\put(6.2,0.15){\line(3,2){2.1}}
\put(6.2,-0.15){\line(3,-2){2.1}}

\put(2.7,-1.5){$2$}
\put(5.7,-1.5){$1$}
\put(9.2,1.8){$3$}
\put(9.2,-2.3){$4$}

\hspace*{9cm}
$ \left( \begin{array}{cccc}
0 & 2 & 2 & 2 \\
3 & 1 & 4 & 4 \\
3 & 4 & 1 & 4 \\
3 & 4 & 4 & 1
\end{array} \right)$

\end{picture}

\vspace*{2cm}

\begin{center} Figure 7 \end{center}

\noindent
The edge $\{1,2\}$, for instance, can be considered to be oriented (in the sense that
its endvertices are situated differently according to the whole graph). 
Algorithm WL applied to this graph yields not only a coloration, but also an 
orientation of the edges, i.e.\ in the colored complete directed graph, which 
represents the cellular algebra, the arcs $(1,2)$ and $(2,1)$ have different colors 
(the colored complete directed graph is given above by its adjacency matrix).
The cellular algebra of $\Gamma$ is of rank $5$. The standard basis consists of
the following matrices:

\[ \left( \begin{array}{cccc}
1 & 0 & 0 & 0 \\
0 & 0 & 0 & 0 \\
0 & 0 & 0 & 0 \\
0 & 0 & 0 & 0
\end{array} \right),
\left( \begin{array}{cccc}
0 & 0 & 0 & 0 \\
0 & 1 & 0 & 0 \\
0 & 0 & 1 & 0 \\
0 & 0 & 0 & 1
\end{array} \right),
\left( \begin{array}{cccc}
0 & 1 & 1 & 1 \\
0 & 0 & 0 & 0 \\
0 & 0 & 0 & 0 \\
0 & 0 & 0 & 0
\end{array} \right),
\left( \begin{array}{cccc}
0 & 0 & 0 & 0 \\
1 & 0 & 0 & 0 \\
1 & 0 & 0 & 0 \\
1 & 0 & 0 & 0
\end{array} \right),
\left( \begin{array}{cccc}
0 & 0 & 0 & 0 \\
0 & 0 & 1 & 1 \\
0 & 1 & 0 & 1 \\
0 & 1 & 1 & 0
\end{array} \right).
\]

Note that for the same graph the algorithm which was described in
\cite{RueR90b} will get as the output the symmetric matrix 
\[
\left(
  \begin{array}{*4{c}}
    0 & 2 & 2 & 2\\
    2 & 1 & 3 & 3\\
    2 & 3 & 1 & 3\\
    2 & 3 & 3 & 1
  \end{array}
\right),
\]
hence two antisymmetric basic matrices of $W(\Gamma)$ will be
merged.

\section{Graph Theoretical Interpretation}

In the previous sections, we already indicated that each cellular
algebra $W$ with basis $A_0,A_1,\ldots,$ $A_{r-1}$ can be represented
by a colored complete directed graph $\Delta=(\Omega,R)$, the graph
whose adjacency matrix is the matrix $A(W)=A(\Delta)=\sum_{i=0}^{r-1}
t_i A_i$ (with indeterminates $t_i$ standing for the colors $i$).
This matrix sometimes is called {\it generic matrix\/} of $W$. For
convenience, the vertices in $\Delta$ have been identified with the
corresponding loops.  In this sense, the basis matrix $A_k$
corresponds to the arc set $R_k$ of color $k$, and $R=\left\{R_0,
  R_1,\ldots,R_{r-1}\right\}$.  The number of colors in $\Delta$ is
equal to the rank of $W$.  More generally, to each linear subspace $S$
with linear basis $A_0,A_1,\ldots,A_{r-1}$ satisfying the properties
(i)-(iv) of a cellular algebra, we can associate a colored complete
directed graph $\Delta=(\Omega,R)$, the graph which belongs to the
generic matrix $\sum_{i=0}^{r-1} t_i A_i$ of $S$.

With this representation, the main idea of the algorithm described above
can be sketched in a more illustrative manner. 
In each iteration of the algorithm, the coloring of the underlying 
complete directed graph is modified by means of Schur-Hadamard 
multiplication and matrix multiplication. These two operations may be interpreted 
as follows. Given generic matrices $X=(x_{uv})$ and $Y=(y_{uv})$ of two 
colored complete directed graphs $\Delta'$ and $\Delta''$ with indeterminates
representing the colors, the Schur-Hadamard product 
$X \circ Y=(x_{uv}y_{uv})$ corresponds to the generic matrix of a new colored 
complete directed graph $\Delta$ where the color of arc $(u,v)$ is the ordered 
mixture of the colors of both arcs in the original graphs. 
In the case of the matrix product $X \cdot Y$, the color of the arc $(u,v)$ in the
new graph $\Delta$ depends on the number and colors of paths of length $2$ starting 
in $u$ and ending in $v$ such that the first step in the path is an arc of $\Delta'$ 
and the second step is an arc of $\Delta''$. The $(u,v)-$entry $\sum_wx_{uw}y_{wv}$ 
of $X \cdot Y$ completely describes the set of these paths.

In the following we present a slightly different procedure leading to
our main goal, namely an efficient computer program 
for graph stabilization. The method is based on the computation of the structure 
constants $p_{ij}^k$ which have been defined and interpreted in Section 4.
To recall the main result, a basis $A_0,A_1,\ldots,A_{r-1}$ of a cellular algebra $W$ 
must fulfill
\[A_i A_j = p_{ij}^0 A_0 + p_{ij}^1 A_1 + \ldots + p_{ij}^{r-1} A_{r-1}\]
for each pair $i,j \in \{0,1,\ldots,r-1\}$.
In the colored graph $\Delta$, this means that each arc $(u,v)$ of a given color $k$ 
is the basis arc of exactly $p_{ij}^k$ triangles with first nonbasis arc of 
color $i$ and second nonbasis arc of color $j$ (a triangle consists of three
not necessarily distinct vertices $u,v,w$ and arcs $(u,v)$, $(u,w)$ and $(w,v)$. 
The arc $(u,v)$ is called the {\it basis arc\/}, the other arcs are the {\it 
nonbasis arcs\/} of the triangle; see Figure 5).

The idea of the algorithm can be described informally as follows.
One iteration includes the round along all arcs of the given graph $\Delta$.
For each arc $(u,v)$ of a fixed color $k$ we count the number of triangles with 
basis arc $(u,v)$ and nonbasis arcs of color $i$ and $j$, respectively, 
$i,j=0,1,\ldots,r-1$ (equivalently, we count the number of paths of length $2$ 
such that the first arc $(u,w)$ is of color $i$ and the second arc $(w,v)$ is of 
color $j$). These numbers should be equal for all arcs. If this is true, then 
these numbers are just the structure constants $p_{ij}^k$. If not, then 
the arc set $R_k$ of color $k$ has to be partitioned into subsets 
$R_{k_0},R_{k_1},\ldots,R_{k_{t-1}}$, each consisting of arcs with the same 
numbers. This step is performed for all colors $k \in \{0,1,\ldots,r-1\}$.
Then the graph $\Delta$ is recolored, i.e.\ we identify color $k_0$ with the 
old color $k$ and introduce the new colors $k_1,\ldots,k_{t-1}$
(in algebraic language, recoloring $\Delta$ means to replace the basis 
matrix $A_k$ by new basis matrices $A_{k_0},A_{k_1},\ldots,A_{k_{t-1}}$).

The next iteration is performed for the recolored graph $\Delta$. If in some 
iteration no new colors are introduced, then the process is stable and we can stop.
In this case, the graph $\Delta$ with the final {\it stable coloring\/} represents 
the required cellular algebra $W$. Here is a more formal description of the algorithm.
\vspace*{0.2cm}

{\bf Algorithm STABIL}

Input: the adjacency matrix $A(\Gamma)=(a_{uv})$ of a colored graph $\Gamma$.

Output: a complete directed graph $\Delta=(\Omega,R)$ with a stable coloring\\
\hspace*{1.4cm} and the structure constants $p_{ij}^k$.

(0) Let $\{0,1,\ldots,s-1\}$ be the set of different entries of $A(\Gamma)$ and\\
\hspace*{0.5cm} $\Delta=(\Omega,R)$ the colored complete directed graph belonging 
to $A(\Gamma)$. 
\\  \hspace*{0.5cm} Determine the arc sets $R_0,R_1,\ldots,R_{s-1}$ of colors 
$0,1,\ldots,s-1$. \\ 
\hspace*{0.5cm} Let $r := s$.

(1) For $k=0,1,\ldots,r-1$ do

\hspace*{1.1cm} For all $(u,v) \in R_k$ do\\
\hspace*{1.7cm} Compute the numbers  $p_{ij}^k$ of triangles with basis arc
$(u,v)$ and \\
\hspace*{1.7cm}  nonbasis arcs of colors $i$ and $j$, respectively,
 $i,j \in \{0,1,\ldots,r-1\}$.

\hspace*{1.1cm} Collect all arcs having the same parameters $p_{ij}^k$, i.e.\
all arcs\\
\hspace*{1.1cm} which belong to the same number of triangles of any colors,\\
\hspace*{1.1cm} and assign them to new sets $R_{k_0},R_{k_1},\ldots,R_{k_{t-1}}$.

\hspace*{1.1cm} Replace $R_k$ by $R_{k_0},R_{k_1},\ldots,R_{k_{t-1}}$,\\
\hspace*{1.1cm} i.e.\ recolor the arcs of color $k$ using the old color
                $k_0=k$ \\
\hspace*{1.1cm} and the new colors $k_1,\ldots,k_{t-1}$.

(2) Let $s$ be the number of colors used to recolor $\Delta$.\vspace*{0.2cm}\\
\hspace*{0.5cm} If $s>r$ then \\
\hspace*{1.1cm} r := s. Goto (1).

(3) STOP.

Let us illustrate this method by a small example, namely we consider
again cuneane, see Figure 3.

Here 
\[A=\left(
  \begin{array}{*8{c}}
    1 & 2 & 3 & 3 & 3 & 3 & 2 & 2\\
    2 & 1 & 2 & 2 & 3 & 3 & 3 & 3\\
    3 & 2 & 1 & 2 & 3 & 3 & 3 & 2\\
    3 & 2 & 2 & 1 & 2 & 3 & 3 & 3\\
    3 & 3 & 3 & 2 & 1 & 2 & 2 & 3\\
    3 & 3 & 3 & 3 & 2 & 1 & 2 & 2\\
    2 & 3 & 3 & 3 & 2 & 2 & 1 & 3\\
    2 & 3 & 2 & 3 & 3 & 2 & 3 & 1
  \end{array}
\right)
\]
is the adjacency matrix of the corresponding colored graph;
\[B=\left(
  \begin{array}{*8{c}}
    1 & 2 & 3 & 4 & 4 & 3 & 2 & 2\\
    2 & 1 & 5 & 5 & 4 & 6 & 4 & 3\\
    3 & 5 & 1 & 5 & 4 & 4 & 6 & 2\\
    4 & 5 & 5 & 1 & 2 & 4 & 4 & 4\\
    4 & 4 & 4 & 2 & 1 & 5 & 5 & 4\\
    3 & 6 & 4 & 4 & 5 & 1 & 5 & 2\\
    2 & 4 & 6 & 4 & 5 & 5 & 1 & 3\\
    2 & 3 & 2 & 4 & 4 & 2 & 3 & 1
  \end{array}
\right)
\]
is the result after the first iteration,
\[
C=\left(
  \begin{array}{*8{c}}
    1 & 2 & 3 & 4 & 4 & 3 & 2 & 5\\
    6 & 7 & 8 & 9 & 10 & 11 & 12 & 13\\
    13 & 8 & 7 & 9 & 10 & 12 & 11 & 6\\
    14 & 15 & 15 & 16 & 17 & 18 & 18 & 14\\
    14 & 18 & 18 & 17 & 16 & 15 & 15 & 14\\
    13 & 11 & 12 & 10 & 9 & 7 & 8 & 6\\
    6 & 12 & 11 & 10 & 9 & 8 & 7 & 13\\
    5 & 3 & 2 & 4 & 4 & 2 & 3 & 1
  \end{array}
\right)
\]
is the result which we get after the second iteration (it in fact
coincides with the final result).



\section{Program Implementation}
The presented algorithm STABIL has been coded in 
C and was tested on a SUN-Sparcstation. The program requires as an
input a file containing the number of colors, the vertex number $n$
and the adjacency matrix $A(\Gamma)$ of an arbitrary graph $\Gamma$,
and provides as an output the number of colors (i.e.\ the rank), the
number of cells in the standard partition, the adjacency matrix of
the cellular algebra $W(\Gamma)$, and (if requested) the structure
constants of $W(\Gamma)$.

In the following we will give some information about the implementation. The adjacency matrix of the graph $\Gamma$ is stored in a 
$n \times n-$matrix. After each iteration of the program, this matrix will contain 
the adjacency matrix of the actual colored complete graph $\Delta$. In the final 
state it contains the adjacency matrix of $W(\Gamma)$. 
 
During any iteration of the program, except the last one, the number of colors 
increases and some arcs of $\Delta$ are recolored. Arising of new colors implies 
a new iteration, while absence of new colors during some iteration gives the sign for
finishing the program. In the latter case a new iteration will not change the 
coloring of $\Delta$.

Any iteration includes the following two imbedded loops: the loop around the 
colors and the loop around the arcs of a fixed color. To handle the loop around 
all arcs of a fixed color, we introduce an additional data structure representing 
the graph $\Delta$, namely lists which store the arcs $(u,v)$ of a given color. 
Each element of a list contains the entries $u$ and $v$, i.e.\ the number of the 
row and column of the arc in the adjacency matrix, and a pointer to the 
information concerning the next arc of the given color.

The following actions with one arc $(u,v)$ of color $k$ form an elementary step 
of the program, except the last one:

(i) Computing the structure constants\footnote{We stress the reader's
  attention that the term ``structure constants'' has a rigorous
  meaning only after the fulfillment of the program. Currently the
  structure constants are just the numbers of triangles with the
  prescribed properties.} $p_{ij}^k$ for $(u,v)$. Since 
$i,j \in \{0,1,\ldots,r-1\}$, i.e.\ $\Delta$ is colored by $r$ colors,  
there are $r^2$ such numbers.\\
(ii) Saving the nonzero structure constants $p_{ij}^k$ for $(u,v)$ as
the sequence of triples $(i,j,p_{ij}^k)$ at the end of a vector MEMORY.\\
(iii) Assigning the color to the arc $(u,v)$. This is done by examining whether
the last sequence in MEMORY is a new sequence (then a new color is introduced),
or the same sequence already appears in MEMORY for some other arc of the 
actual color $k$.

Let us consider these actions in some more detail.

{\it Computation of the structure constants\/}. As mentioned above, each arc
is characterized by a set of $r^2$ numbers, called structure constants. The
geometrical meaning of these numbers has been described, too. In order to
calculate the numbers $p_{ij}^k$ for a given arc $(u,v)$, we examine all
triangles with basis arc $(u,v)$ (see Figure 5). Note that there are $n$ such 
triangles. If for some vertex $w$ the arc $(u,w)$ has color $i$ and the arc $(w,v)$
has color $j$, then we increase $p_{ij}^k$ by $1$ (initially all $p_{ij}^k$ 
are set equal to $0$).  Thus the number of nonzero $p_{ij}^k$ does not
exceed $n$.

If we compute simultaneously all $r^2$ structure constants for an arc and store 
them in a straightforward way using an $r \times r-$matrix, then we will soon 
get storage overflow, since the existence of a matrix with $r^2$ elements in the 
program is impossible already for comparatively small $n$. However, since
there are at most $n$ nonzero structure constants for each arc, we can use 
instead of a matrix CONST with $r^2$ entries a data structure whose size is 
proportional to $n$. This data structure consists of lists whose elements contain 
the information $i$, $p_{ij}^k$ and a pointer to the next nonzero structure 
constant in the column $j$ of the matrix CONST. Additionally we need a vector with 
$r$ elements to save the pointers to the columns 
(the described technique, called {\it hashing\/},
is a well known tool to treat efficiently set manipulation problems; see
\cite{AhoHU74}).

{\it Saving of the structure constants for arcs of a given color\/}.
In view of the forthcoming manipulations, it is more adequate to save the structure 
constants as a sequence of triples $i,j,p_{ij}^k$. The maximal number of such 
triples in the sequence of a given arc is $n$. During the pass along the arcs of 
a given color $k$, we have to save all sequences of structure constants belonging to
the arcs of this color. Therefore the length of the vector MEMORY should be 
about $3 n^3$ (note that the number of arcs of color $k$ is restricted by $n^2$).

{\it The search along MEMORY and recoloring of an arc\/}.
The nonzero structure constants for a given arc $(u,v)$ of color $k$ are saved 
at the end of the vector MEMORY. Then we should examine whether we have a new
sequence or whether an identical sequence has been saved in MEMORY before. 
We use the partial ordering of the sequences by their lengths. So, if we search 
for an identical sequence, in the case of equal lengths we come to an element 
by element comparison, in the case of nonequal lengths we pass to the next sequence.

As mentioned before, the main problem in the program is to find memory for
the structure constants. For graphs with a comparatively large number of 
vertices, it is impossible to save all these constants for a given color $k$.
Therefore we decided to save only the nonzero numbers as triples 
$i,j,p_{ij}^k$ in 
the vector MEMORY. Now, some graphs may initially produce a very large number
of different sequences $i,j,p_{ij}^k$ and the corresponding massive 
storage requirements may exceed the available memory capacity of small computers.
Therefore we start the program with a simple {\it preprocessing procedure\/} 
which aims to increase in advance the number of colors, before starting the main
program.

In the first step of this preprocessing we classify the vertices of the 
(colored) graph $\Gamma$ in the following way. Two vertices are put into the 
same cell if and only if they are incident to the same number of edges of each color. 
Then we can recolor the edges according to the new coloring of the vertices.
If two edges, which initially have the same color, connect two different colored
pairs of vertices, then the edges are assigned different colors, too.

Some additional practically important details about the current
version of the program implementation of algorithm STABIL may be
found in Section 9.

\section{Estimation of the Complexity}

The question of the theoretical complexity\footnote[1]{Readers who are not 
familiar with complexity considerations of algorithms are referred to the 
standard book \cite{AhoHU74}.} of the WL--stabilization was not 
considered for a long time. Weisfeiler and Leman only stated that the complexity 
is polynomial in the vertex number $n$ of the graph, without giving any explicit 
time bound. A first attempt for an estimation was done by S.\ Friedland. In 
\cite{Fri89} he pointed out that the required time is restricted by $O(n^{10})$.
I.N.\ Ponomarenko \cite{Pon93a} improved the time bound to $O(n^5 \, log \,n)$.
Very recently, L.\ Babel showed in \cite{Bab95} that the algorithm can be implemented 
to run in time $O(n^3 \, log \,n)$.

{\bf 8.1.}
Before we are going to determine the worst case complexity of our implementation, 
let us once more stress that we are dealing with two different problems,
depending on the point of view. The first problem 
is to find the standard basis of the cellular algebra $W(\Gamma)$ belonging to
some graph $\Gamma$. The second problem, which many times appears in
framework of algebraic combinatorics, is to compute the structure constants $p_{ij}^k$ and
a stable coloring of the complete directed graph $\Delta$ representing the 
cellular algebra $W(\Gamma)$. Our implementation solves the second problem 
whereas the implementation of the above-mentioned algorithm by Babel
solves the first problem.

Let us now analyze the complexity of algorithm STABIL. We only have to examine step (1),
the initializing step (0) and step (2) obviously can be performed in time $O(n^2)$.
We have seen in the previous section that an elementary step of (1) consists of 
three actions which are performed for each of the $n^2$ arcs of the graph.

In action (i) the $n$ triangles with basis arc $(u,v)$ have to be found. This can 
be done in time $O(n)$ by inspecting the $u$th row and the $v$th column of the 
adjacency matrix $A=(a_{uv})$ of the graph. Note that for each vertex $w$ there 
is one triangle with basis arc $(u,v)$, the entries $a_{uw}$ and $a_{wv}$ of $A$ 
are the colors of the nonbasis arcs $(u,w)$ and $(w,v)$.
For each triangle the value of some parameter $p_{ij}^k$ is actualized. 
To be more precise, let $a_{uv}=k$. Then the value of $p_{ij}^k$ must be increased 
by 1 if $a_{uw}=i$ and $a_{wv}=j$.
In order to find the actual value of $p_{ij}^k$ we have to pass through the 
corresponding list in the data structure (containing the information $i, p_{ij}^k$).
Since the length of this list is at most $n$, this requires time 
$O(n)$ for each triangle. Thus, action (i) requires total time $O(n^2)$ for one arc.
Since there are at most $n$ nonzero structure constants for each arc, action (ii),
namely saving the sequence of structure constants in the vector MEMORY, 
can be executed in time $O(n)$. The most time consuming part is action (iii).
In order to compare the sequence of triples $i,j,p_{ij}^k$ for one arc $(u,v)$
with all such sequences already stored in MEMORY, we eventually have to pass 
through the whole vector MEMORY. Since this vector may be of length $3n^3$,
this requires time $O(n^3)$. (We stress that the current implementation does
not use any storage/search technique, see also 8.5).

This shows that the complexity of one elementary step is $O(n^3)$. Since $n^2$
arcs are treated, one iteration of step (1) requires time $O(n^5)$.
Now it remains to give a bound on the number of iterations.
If only one color is added during each iteration, then there are $n^2$ iterations.
Altogether, this results in a worst case time bound of $O(n^7)$.

\noindent{\bf Remarks}
\begin{enumerate}
\item Similar reasonings were done less carefully in \cite{ChuKP92}
  and thus resulted in the evaluation $O(n^8)$.
\item In fact, the number of iterations in the WL-stabilization is
  less than $n^2$ (\cite{Ade95}), however such an opportunity to
  improve the evaluation will not be used in this paper.
\item In the worst case, if $W(\Gamma)$ coincides with the full matrix
  algebra of order $n$, there are $n^2$ basis matrices and therefore
  $n^6$ structure constants (most of them are zero). This information
  may help the reader to realize the difference between the statements
  of Problems 1 and 2.
\end{enumerate}

The crucial point in the implementation concerning both running time
and space requirement is the vector MEMORY. In order to make the
program applicable also for relatively small computers, it is
favourable to restrict the length of MEMORY to $O(n^2)$\footnote[1]{In
  the actual version of the program, we defined MEMORY to be of length
  $3n^2$, see also additional remarks in next section}. With this
modification, it may be impossible to store all the different
sequences of structure constants. In that case, new sequences are not
saved and all the corresponding arcs are assigned the same color. This
color will be split during the next iteration of the program (note
that this procedure may increase the number of iterations which are
needed to obtain the stable coloring).

{\bf 8.2.}  As mentioned above, the implementation of the algorithm
presented in \cite{Bab95} has a considerably lower worst case
complexity. We will very briefly indicate the main ideas of that
implementation. Compared to algorithm STABIL, there are two main
modifications, one decreases the number of triangles which are
examined in one iteration, the other involves some sophisticated
sorting techniques.

In each iteration of algorithm STABIL, all $n^3$ triangles of the colored graph 
$\Delta$ are examined. However, it is not really necessary to inspect the whole 
set of triangles. One can restrict to a certain subset. 
Roughly sketched, the procedure is the following.
Let a colored complete directed graph $\Delta$ be given. 
In each iteration some arcs of the graph will keep their colors, others are assigned
new colors (which have not been used in the previous iteration). More concrete, the 
arc set $R_k$ is split into subsets $R_{k_0},R_{k_1},\ldots,R_{k_{t-1}}$, where one 
of the colors $k_0,k_1,\ldots,k_{t-1}$ is equal to $k$ and the others are new colors.
The basic idea is to inspect only those triangles which contain at least one arc of
a new color. Denote by $T_M$ the set of these triangles. Further let $R_M$ 
denote the set of arcs which are basis arcs of triangles from $T_M$.
Now, step (1) of the algorithm STABIL is modified in the following way.
Each arc $(u,v) \in R_M$ is the basis arc of some (in general less than $n$) triangles 
from $T_M$. For $(u,v) \in R_M$ list the colors of the nonbasis arcs of these triangles.
Now, each arc from $R_M$ is associated a multiset of some ordered pairs $(i,j)$.
Collect arcs with equal multisets and assign them the same color, i.e.\ replace
each arc set $R_k$ by suitable subsets $R_{k_0},R_{k_1},\ldots,R_{k_{t-1}}$. One
of these subsets, say $R_{k_0}$ (which may be empty), consists of all arcs from $R_k$ 
which do not belong to $R_M$. 
The procedure stops if no new colors are generated.

It is not yet designated which one of the colors $k_0,k_1,\ldots,k_{t-1}$ is equal to
the old color $k$ and which ones are new colors. The effort for each iteration is kept
low if $T_M$ contains only a small number of triangles. Therefore it is favourable
to identify $k$ with that color $k_p$ where $R_{k_p}$ contains the largest number
of arcs. 
It is not hard to check the correctness of this method (for details see
\cite{Bab95}). The worst case complexity is determined as follows.

Let $\tau_h$ denote the cardinality of $T_M$ in the $h$th iteration.
Then obviously $|R_M| \leq \tau_h$.
Therefore, multisets of at most $\tau_h$ arcs $(u,v)$ have to be computed.
These multisets are stored as lists $S(u,v)$ and are obtained as follows.
Order the $\tau_h$ triangles from $T_M$ lexicographically according to the colors 
$(i,j)$ of the
nonbasis arcs (a pair $(i,j)$ appears before $(i',j')$ if and only if $i<i'$ or 
$i=i'$ and $j \leq j'$). Note that the colors are in the range $\{0,1,\ldots,n^2-1\}$
(since the graph has $n^2$ arcs, not more than $n^2$ colors can occur).
It is well known that lexicographical ordering of $\tau_h$ pairs of integers from 
$\{0,1,\ldots,n^2-1\}$ can be done in time $O(\tau_h+n)$ using the sorting routine 
{\it bucket sort\/}.
Now, the lists $S(u,v)$ are obtained by passing through the ordered list of triangles
and assigning the actual triangle to its basis arc $(u,v)$, i.e.\ the colors of the
nonbasis arcs are inserted at the end of $S(u,v)$. Obviously, this requires time 
$O(\tau_h)$. Note that the pairs of colors in the lists $S(u,v)$ appear now in 
lexicographical order.

To identify the different multisets we have to order the lists $S(u,v)$, 
$(u,v) \in R_M$, lexicographically. Since the total length of all lists is $\tau_h$
and the entries are pairs of numbers from $\{0,1,\ldots,n^2-1\}$, this again can be 
done with bucket sort in time $O(\tau_h+n)$.
Now the arc sets $R_k$ are split in the obvious way by passing through the 
ordered list of multisets $S(u,v)$. This, as well as finding the subsets $R_{k_p}$ of 
largest cardinality, requires time $O(\tau_h)$.

So far we have seen that the complexity for the $h$th iteration of this method
is $O(\tau_h+n)$. It remains to compute the total complexity for all iterations 
(trivially, the number of iterations is restricted by the maximal number $n^2$ of 
colors).
 
Since $R_{k_p}$ has been chosen to be that subset of $R_k$ with largest cardinality,
each of the other subsets (which contain the arcs with new colors) has at most half 
the size of $R_k$.
Therefore, each time a certain triangle is inspected, at least one arc set which 
shares an arc with this triangle is at most half as large than before.
As a consequence, each of the $n^3$ triangles is examined not more than 
$3\,log\,n^2 = 6\,log\,n$ times. This shows that $\sum_h \tau_h \leq 6\,n^3\,log\,n$.
Finally we obtain a worst case time complexity of $O(n^3\,log\,n)$.

{\bf 8.3.}
These ideas have been realized in a computer program by L.\ Babel, S.\ Baumann
and M.\ L\"udecke. The program is termed STABCOL, due to the fact that the coloring 
of the complete directed graph is modified in each iteration until the process 
is stable, i.e.\ until a STABle COLoring is obtained. It is coded in programming 
language C and was also tested on the same SUN-Sparcstation. Just as for the program STABIL, the 
input is a file containing the number of colors, the vertex number and the 
adjacency matrix of a graph $\Gamma$, the output contains the number of 
colors, the number of cells and the adjacency matrix of the cellular algebra 
which is generated by $\Gamma$.

Contrary to STABIL, the program STABCOL does not work with vectors of predefined lengths
but uses more sophisticated data structures. The set $T_M$ of triangles and the 
multisets $S(u,v)$ are stored in lists which are linked by pointers and which are 
of variable length. In this way, waste of memory space is avoided. Furthermore,
memory space which is no longer needed is set free immediately.

The complexity analysis of STABIL shows that most of its time is spent in order to
compare a new sequence of numbers with old sequences. Since this is done in the
obvious way by passing through the entire vector of sequences, it requires time 
proportional to the length of the vector. This somewhat time-consuming 
procedure is 
avoided in STABCOL by means of very special sorting techniques. 
These techniques and the more complicated data structures make the implementation
much more ambitious. However, we do not have enough space to go into details here.
The interested reader may consult the program description \cite{BabBLT97}.

{\bf 8.4.}
At first glance, a comparison of the theoretical complexities indicates that
the implementation of \cite{Bab95} should be preferred. However, it turns out that 
this implementation, although theoretically very fast, is rather slow in practice
and applicable only for relatively small graphs, whereas our implementation, although 
inferior with respect to the worst case bound, is very fast in practice
and is able to 
handle very large graphs (the practical behaviour of both program implementations is 
documented in the next section).

Here we are confronted with a situation which seems to be strange but which rather 
frequently occurs on the construction of algorithms. There are two algorithms or 
two implementations of an algorithm solving the same problem, one of them 
theoretically fast (i.e.\ with a good worst case complexity) but practically slow, 
the other one practically fast in spite of a relatively bad worst case complexity.

There are two main reasons for this paradox. First, and perhaps most
important is that the sign ``O'' in the evaluation means in fact the
existence of some constant as a multiplier with the monomial depending on
$n$. The actual value of this constant depends on many factors, in
particular on the ``complexity'' of the data structures. In our case
the multiplier for the implementation of STABIL is essentially smaller
than the one for STABCOL. Thus the advantages of the theoretically
faster algorithm cannot be felt on comparably small graphs. Second,
one algorithm has been constructed from a purely theoretical point of
view with the aim to obtain a worst case complexity as good as
possible. No practical considerations are taken into account such as
simplicity of data structures, easy way of implementing, small space
requirements, etc. In particular, the running time of the algorithm in
the mean (the mean taken over a large representative selection of
practically relevant examples) is not considered.  This average
behaviour, however, is much more important for practitioners than the
worst case behaviour, which often occurs only for pathological
examples.  The second algorithm is constructed from that practical
point of view. It aims to solve the ``real world'' problems very fast,
without paying attention to its theoretical complexity.

{\bf 8.5.}
It is worthwhile to stress that the worst case time bound $O(n^{7})$ (see 8.1)
is a very rough upper bound.

In fact we see a number of opportunities to diminish this bound essentially.
One of them was mentioned  in Remark (ii) in Subsection 8.1.

Also (as a result of restricting the length of vector MEMORY to $O(n^{2})$) the
complexity of one elementary step
in the current version of the program is actually reduced from $O(n^{3})$ to
$O(n^{2})$. In spite of the fact that the total number of steps may slightly
grow, here we have one more standby to reduce the upper bound.

A more careful analysis of the applied technique of hashing together with a more
clever organization of storage (via the use, e.g., of balanced binary trees) may
essentially decrease the number of comparisons when we operate with the vector MEMORY.

However we do not use these and other possible options in the current
preliminary version of our report. In contrast to \cite{BabBLT97}, our report
is oriented towards those practical users of STABIL, for whom the theoretical
question of the evaluation of the efficiency does not play a crucial role.

Nevertheless, we intend to return to the consideration of this question in the future.

\section{Testing the Program}

The presented implementation STABIL has been tested on a large number of structures.
All computations were done on a SUN-Sparcstation 10 with 128MB RAM.
On this machine, the program is able to handle graphs with up to 2000 vertices.
To demonstrate the capability of the program we considered acyclic 
compounds and compounds containing multiple bonds or heteroatoms as depicted in Figure 8.
These structures also appeared as illustrations 
in papers by other authors (see e.g.\ \cite{RBW80}, \cite{RueR90b}).
The results are summarized in Table I. Besides the running time of the program,
the number of cells in the standard partition and the number of
colors in the stable coloring (i.e.\ the number of equivalence classes 
of atoms and ordered pairs of atoms) are stated. In order to make evident the 
practical efficiency of our program, we also state the running time of the 
program implementation STABCOL. Note that STABCOL requires space proportional to $n^3$,  
therefore it can handle graphs with not more than 150 vertices.


\noindent
\hspace*{3.5cm}
{\bf Table I}. Results for the structures in Figure 8

\vspace*{0.4cm}

\hspace*{2.5cm}
\begin{tabular}{cccccc} \hline
structure &   & number   & number    & \multicolumn{2}{c}{CPU time [seconds]} \\
(graph)   & n & of cells & of colors & STABIL & STABCOL \\ \hline
1 &  12 & 2  & 16  & 0.04 & 0.05 \\
2 &  8  & 3  & 18  & 0.03 & 0.01 \\
3 &  12 & 7  & 66  & 0.05 & 0.07 \\
4 &  18 & 6  & 86  & 0.08 & 0.27 \\
5 &  12 & 12  & 144 & 0.05 & 0.08 \\   
6 &  20 & 1  & 6   & 0.06 & 0.19 \\
7 &  12 & 3  & 27  & 0.04 & 0.07 \\
8 &  18 & 8  & 102 & 0.07 & 0.25 \\
9 &  17 & 17  & 289 & 0.06 & 0.21 \\
10 & 20 & 3  & 30  & 0.07 & 0.32 \\
11 & 18 & 8  & 83 & 0.07 & 0.21 \\
12 & 20 & 10  & 119 & 0.09 & 0.30 \\
13 & 28 & 27  & 730 & 0.20 & 1.00 \\
14 & 10 & 5  & 34  & 0.04 & 0.03 \\
15 & 11 & 3  & 22  & 0.04 & 0.05 \\
16 & 22 & 8 & 146 & 0.13 & 0.43 \\
17 & 20 & 6  & 76  & 0.08 & 0.31 \\
18 & 16 & 10  & 114 & 0.06 & 0.16 \\
19 & 8  & 2  & 8  & 0.03 & 0.01 \\        
20 & 22 & 16  & 292 & 0.14 & 0.45 \\
21 & 25 & 5  & 31  & 0.08 & 0.34 \\
22 & 8  & 8  & 64  & 0.03 & 0.01 \\ \hline
\end{tabular}

\includegraphics{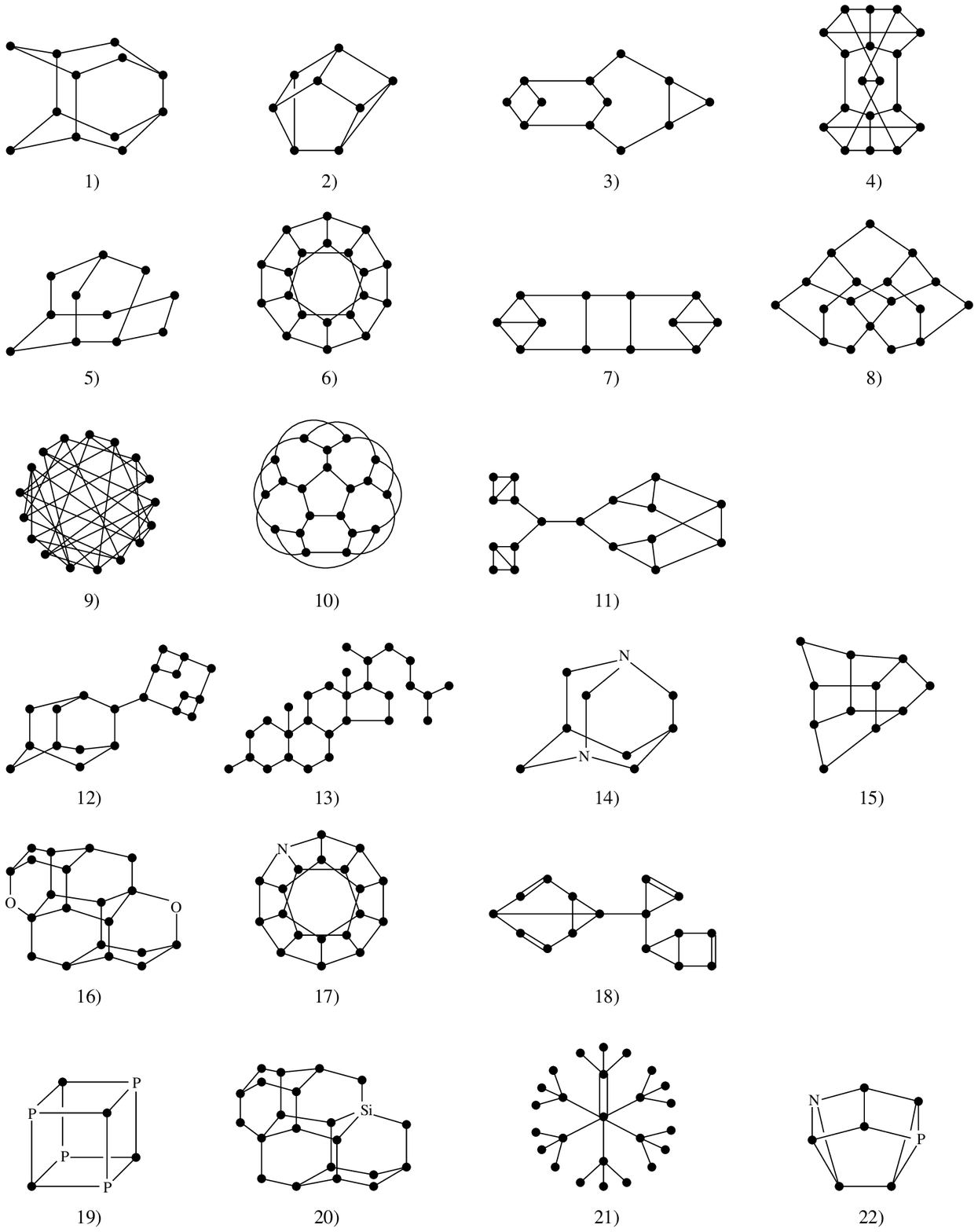} 

\vspace*{0.5cm}

\begin{center} Figure 8 \end{center}

\newpage

The program has also been tested on first members of three infinite families of 
graphs where the automorphism groups and the numbers of orbits on the vertices
and ordered pairs of vertices are known.
We give a description of these families. The results are shown in Tables II--IV.

{\bf Benzene stacks}.
We denote by $P_k$ the graph from this family consisting of $n=6k$ vertices. 
The vertices of $P_k$ form $k$ stages (strata), each stage (stratum) is a cycle 
of $6$ vertices. Besides the edges in these cycles there are edges between
stages. The graphs $P_k$, $k=2,3,4$, are depicted in Figure 9.
A formal description of the graphs $P_k$ is the following.

Let $L=\{a,b,c,d,e,f\}$, $K=\{1,2,\ldots,k\}$, $x_i=(x,i)$ for $x \in L, i \in K$.\\
Then $P_k=(\Omega(P_k),E(P_k))$, where $\Omega(P_k)=L \times K$,
$E(P_k)=\bigcup_{i=1}^k R_i \cup \bigcup_{j=1}^{k-1} Q_j$, and\\
$R_i=\{\{a_i,b_i\},\{b_i,c_i\},\{c_i,d_i\},\{d_i,e_i\},\{e_i,f_i\},\{f_i,a_i\}\}$,\\
$Q_j=\left\{\begin{array}{ccl}
    \{\{a_j,a_{j+1}\},\{c_j,c_{j+1}\},\{e_j,e_{j+1}\}\} & : & j=2l-1 \\
    \{\{b_j,b_{j+1}\},\{d_j,d_{j+1}\},\{f_j,f_{j+1}\}\} & : & j=2l.
  \end{array} \right.$

It is known from \cite{KLP89} and from \cite{KliLPZ92} that
the automorphism group of $P_k$ is isomorphic to $S_3 \times
S_2$. $Aut(P_k)$ has $k$ orbits on the set $\Omega(P_k)$ and
$4k^2$ orbits on the set $\Omega(P_k) \times \Omega(P_k)$.

\vspace*{0.6cm}

\noindent

\begin{center} {\bf Table II}. Results for benzene stacks

\vspace*{0.4cm}

\begin{tabular}{ccccc} \hline
   & number   & number    & \multicolumn{2}{c}{CPU time [seconds]} \\
 n & of cells & of colors & STABIL & STABCOL \\ \hline
6   &  1 &    4 &   0.03 &   0.01 \\
12  &  2 &   16 &   0.03 &   0.08 \\
18  &  3 &   36 &   0.08 &   0.32 \\
24  &  4 &   64 &   0.15 &   0.88 \\
30  &  5 &  100 &   0.28 &   2.01 \\
36  &  6 &  144 &   0.44 &   4.83 \\
42  &  7 &  196 &   0.79 &  10.60 \\
48  &  8 &  256 &   1.14 &  14.90 \\
54  &  9 &  324 &   1.71 &  21.87 \\
60  & 10 &  400 &   2.48 &  24.73 \\
66  & 11 &  484 &   3.39 &  32.77 \\
72  & 12 &  576 &   4.60 &  58.35 \\ 
78  & 13 &  676 &   6.64 &  75.92 \\
102 & 17 & 1156 &  16.21 & 177.82 \\
126 & 21 & 1764 &  35.45 & 362.34 \\
150 & 25 & 2500 &  65.74 & -- \\
174 & 29 & 3364 & 117.89 & -- \\
198 & 33 & 4356 & 190.83 & -- \\ \hline
\end{tabular}
\end{center}

\vspace*{1.3cm} 
\setlength{\unitlength}{0.00083333in}
\begingroup\makeatletter\ifx\SetFigFont\undefined%
\gdef\SetFigFont#1#2#3#4#5{%
  \reset@font\fontsize{#1}{#2pt}%
  \fontfamily{#3}\fontseries{#4}\fontshape{#5}%
  \selectfont}%
\fi\endgroup%
{\renewcommand{\dashlinestretch}{30}
\begin{picture}(6865,2710)(0,-10)
\put(3508,720){\blacken\ellipse{82}{82}}
\put(3508,720){\ellipse{82}{82}}
\put(2665,1205){\blacken\ellipse{82}{82}}
\put(2665,1205){\ellipse{82}{82}}
\put(2935,1354){\blacken\ellipse{82}{82}}
\put(2935,1354){\ellipse{82}{82}}
\put(3505,1023){\blacken\ellipse{82}{82}}
\put(3505,1023){\ellipse{82}{82}}
\put(4077,1354){\blacken\ellipse{82}{82}}
\put(4077,1354){\ellipse{82}{82}}
\put(4344,1205){\blacken\ellipse{82}{82}}
\put(4344,1205){\ellipse{82}{82}}
\put(2665,2172){\blacken\ellipse{82}{82}}
\put(2665,2172){\ellipse{82}{82}}
\put(2935,2015){\blacken\ellipse{82}{82}}
\put(2935,2015){\ellipse{82}{82}}
\put(4086,2020){\blacken\ellipse{82}{82}}
\put(4086,2020){\ellipse{82}{82}}
\put(3505,2347){\blacken\ellipse{82}{82}}
\put(3505,2347){\ellipse{82}{82}}
\put(4348,2166){\blacken\ellipse{82}{82}}
\put(4348,2166){\ellipse{82}{82}}
\put(3508,2044){\blacken\ellipse{82}{82}}
\put(3508,2044){\ellipse{82}{82}}
\put(3819,1871){\blacken\ellipse{82}{82}}
\put(3819,1871){\ellipse{82}{82}}
\put(3819,1504){\blacken\ellipse{82}{82}}
\put(3819,1504){\ellipse{82}{82}}
\put(3508,1326){\blacken\ellipse{82}{82}}
\put(3508,1326){\ellipse{82}{82}}
\put(3190,1500){\blacken\ellipse{82}{82}}
\put(3190,1500){\ellipse{82}{82}}
\put(3190,1868){\blacken\ellipse{82}{82}}
\put(3190,1868){\ellipse{82}{82}}
\put(3505,2646){\blacken\ellipse{82}{82}}
\put(3505,2646){\ellipse{82}{82}}
\path(3506,717)(4345,1201)(4345,2169)
        (3506,2653)(2668,2169)(2668,1201)(3506,717)
\path(3506,1322)(3820,1504)(3820,1867)
        (3506,2048)(3192,1867)(3192,1504)(3506,1322)
\path(3506,2653)(3506,2351)
\path(3506,1020)(4082,1352)(4082,2019)
        (3506,2351)(2930,2018)(2930,1352)(3506,1020)
\path(4343,1203)(4081,1354)
\path(2667,1203)(2928,1354)
\path(3819,1869)(4082,2021)
\path(3506,1322)(3506,1020)
\path(3191,1869)(2928,2019)
\put(888,692){\blacken\ellipse{82}{82}}
\put(888,692){\ellipse{82}{82}}
\put(1726,1179){\blacken\ellipse{82}{82}}
\put(1726,1179){\ellipse{82}{82}}
\put(1726,2150){\blacken\ellipse{82}{82}}
\put(1726,2150){\ellipse{82}{82}}
\put(888,2627){\blacken\ellipse{82}{82}}
\put(888,2627){\ellipse{82}{82}}
\put(49,2148){\blacken\ellipse{82}{82}}
\put(49,2148){\ellipse{82}{82}}
\put(49,1183){\blacken\ellipse{82}{82}}
\put(49,1183){\ellipse{82}{82}}
\put(577,1479){\blacken\ellipse{82}{82}}
\put(577,1479){\ellipse{82}{82}}
\put(888,1298){\blacken\ellipse{82}{82}}
\put(888,1298){\ellipse{82}{82}}
\put(1202,1479){\blacken\ellipse{82}{82}}
\put(1202,1479){\ellipse{82}{82}}
\put(1202,1844){\blacken\ellipse{82}{82}}
\put(1202,1844){\ellipse{82}{82}}
\put(888,2025){\blacken\ellipse{82}{82}}
\put(888,2025){\ellipse{82}{82}}
\put(577,1844){\blacken\ellipse{82}{82}}
\put(577,1844){\ellipse{82}{82}}
\put(5972,2629){\blacken\ellipse{82}{82}}
\put(5972,2629){\ellipse{82}{82}}
\put(6816,2149){\blacken\ellipse{82}{82}}
\put(6816,2149){\ellipse{82}{82}}
\put(6812,1180){\blacken\ellipse{82}{82}}
\put(6812,1180){\ellipse{82}{82}}
\put(5976,692){\blacken\ellipse{82}{82}}
\put(5976,692){\ellipse{82}{82}}
\put(5137,1180){\blacken\ellipse{82}{82}}
\put(5137,1180){\ellipse{82}{82}}
\put(5129,2145){\blacken\ellipse{82}{82}}
\put(5129,2145){\ellipse{82}{82}}
\put(5972,2322){\blacken\ellipse{82}{82}}
\put(5972,2322){\ellipse{82}{82}}
\put(6549,1995){\blacken\ellipse{82}{82}}
\put(6549,1995){\ellipse{82}{82}}
\put(6549,1334){\blacken\ellipse{82}{82}}
\put(6549,1334){\ellipse{82}{82}}
\put(5976,1010){\blacken\ellipse{82}{82}}
\put(5976,1010){\ellipse{82}{82}}
\put(5396,1329){\blacken\ellipse{82}{82}}
\put(5396,1329){\ellipse{82}{82}}
\put(5400,1995){\blacken\ellipse{82}{82}}
\put(5400,1995){\ellipse{82}{82}}
\put(5972,2023){\blacken\ellipse{82}{82}}
\put(5972,2023){\ellipse{82}{82}}
\put(6287,1846){\blacken\ellipse{82}{82}}
\put(6287,1846){\ellipse{82}{82}}
\put(6291,1479){\blacken\ellipse{82}{82}}
\put(6291,1479){\ellipse{82}{82}}
\put(5972,1305){\blacken\ellipse{82}{82}}
\put(5972,1305){\ellipse{82}{82}}
\put(5661,1483){\blacken\ellipse{82}{82}}
\put(5661,1483){\ellipse{82}{82}}
\put(5657,1842){\blacken\ellipse{82}{82}}
\put(5657,1842){\ellipse{82}{82}}
\put(5815,1753){\blacken\ellipse{82}{82}}
\put(5815,1753){\ellipse{82}{82}}
\put(6126,1753){\blacken\ellipse{82}{82}}
\put(6126,1753){\ellipse{82}{82}}
\put(6126,1572){\blacken\ellipse{82}{82}}
\put(6126,1572){\ellipse{82}{82}}
\put(5972,1491){\blacken\ellipse{82}{82}}
\put(5972,1491){\ellipse{82}{82}}
\put(5815,1572){\blacken\ellipse{82}{82}}
\put(5815,1572){\ellipse{82}{82}}
\put(5976,1842){\blacken\ellipse{82}{82}}
\put(5976,1842){\ellipse{82}{82}}
\path(889,697)(1727,1180)(1727,2149)
        (889,2632)(50,2149)(50,1180)(889,697)
\path(889,1302)(1204,1483)(1204,1846)
        (889,2027)(574,1846)(574,1483)(889,1302)
\path(5973,697)(6812,1180)(6812,2149)
        (5973,2632)(5135,2149)(5135,1180)(5973,697)
\path(5973,999)(6549,1331)(6549,1998)
        (5973,2330)(5397,1997)(5397,1331)(5973,999)
\path(5973,1302)(6288,1483)(6288,1846)
        (5973,2027)(5658,1846)(5658,1483)(5973,1302)
\path(5973,1483)(6131,1573)(6131,1756)
        (5973,1846)(5816,1755)(5816,1573)(5973,1483)
\path(889,2632)(889,2027)
\path(1202,1485)(1726,1180)
\path(573,1484)(48,1181)
\path(5973,2632)(5973,2330)
\path(6811,1182)(6548,1333)
\path(5134,1182)(5396,1333)
\path(6288,1848)(6548,1999)
\path(5973,1302)(5973,999)
\path(5657,1847)(5396,1998)
\path(5973,2027)(5973,1846)
\path(6286,1485)(6130,1575)
\path(5657,1485)(5815,1575)
\put(889,31){\makebox(0,0)[lb]{\smash{{{\SetFigFont{9}{10.8}{\rmdefault}{\mddefault}{\updefault}$P_2$}}}}}
\put(3491,31){\makebox(0,0)[lb]{\smash{{{\SetFigFont{9}{10.8}{\rmdefault}{\mddefault}{\updefault}$P_3$}}}}}
\put(5973,31){\makebox(0,0)[lb]{\smash{{{\SetFigFont{9}{10.8}{\rmdefault}{\mddefault}{\updefault}$P_4$}}}}}
\end{picture}
}

\vspace*{0.5cm}

\begin{center} Figure 9 \end{center}

{\bf M{\"o}bius ladders}. 
We denote by $M_k=(\Omega(M_k),E(M_k))$ the graph 
with the set of $n=2k$ vertices
$\Omega(M_k)=\{a_1,\ldots,a_k,a_{k+1},\ldots,a_{2k}\}$ 
and the set of edges
$E(M_k)=\{\{a_i,a_j\} \, | \; j-i=x \, ($mod $2k), x \in \{1,k,2k-1\}
\}$. For example, the graph $M_5$ is depicted in Figure 10, which may
serve as an explanation of the name.

The symmetry of the graphs $M_k$ has been investigated in \cite{KKZ90}, \cite{Sim86}, 
\cite{WSH88}, \cite{FarKM94} and \cite{KliRRT95}.
It was proved in \cite{KKZ90} that, for $k>3$, the automorphism group of $M_k$ is 
isomorphic to the dihedral group $D_{2k}$. This group has one orbit on the set 
$\Omega(M_k)$ and $k+1$ orbits on the set $\Omega(M_k) \times \Omega(M_k)$.

\vspace*{0.6cm}

\noindent

\begin{center} {\bf Table III}. Results for M{\"o}bius ladders

\vspace*{0.4cm}

\begin{tabular}{ccccc} \hline
   & number   & number    & \multicolumn{2}{c}{CPU time [seconds]} \\
 n & of cells & of colors & STABIL & STABCOL \\ \hline
6   & 1 & 3   &  0.03 &   0.01 \\
12  & 1 & 7   &  0.05 &   0.08 \\
18  & 1 & 10  &  0.07 &   0.28 \\
24  & 1 & 13  &  0.11 &   0.80 \\
30  & 1 & 16  &  0.21 &   1.68 \\
36  & 1 & 19  &  0.35 &   3.99 \\
42  & 1 & 22  &  0.57 &   7.62 \\
48  & 1 & 25  &  0.87 &  12.34 \\
54  & 1 & 28  &  1.27 &  17.74 \\
60  & 1 & 31  &  1.76 &  24.34 \\
66  & 1 & 34  &  2.39 &  32.54 \\
72  & 1 & 37  &  3.21 &  46.22 \\
80  & 1 & 41  &  4.57 &  66.16 \\
100 & 1 & 51  &  9.57 & 133.30 \\
120 & 1 & 61  & 17.33 & 232.40 \\
140 & 1 & 71  & 29.03 & 388.04 \\ 
160 & 1 & 81  & 47.03 &  --    \\
180 & 1 & 91  & 69.87 &  --    \\
200 & 1 & 101 & 95.85 &  --    \\ \hline
\end{tabular}
\end{center}

\begin{center}
\setlength{\unitlength}{0.00083300in}%
\begingroup\makeatletter\ifx\SetFigFont\undefined%
\gdef\SetFigFont#1#2#3#4#5{%
  \reset@font\fontsize{#1}{#2pt}%
  \fontfamily{#3}\fontseries{#4}\fontshape{#5}%
  \selectfont}%
\fi\endgroup%
\begin{picture}(2850,3045)(4801,-5500)
\thicklines
\put(6001,-2761){\circle*{150}}
\put(6901,-2761){\circle*{150}}
\put(7501,-3661){\circle*{150}}
\put(7501,-4261){\circle*{150}}
\put(6001,-5161){\circle*{150}}
\put(6001,-4561){\circle*{150}}
\put(5101,-4261){\circle*{150}}
\put(5101,-3661){\circle*{150}}
\put(6001,-3361){\circle*{150}}
\put(6901,-3361){\circle*{150}}
\put(6001,-2761){\line( 1, 0){900}}
\put(6901,-2761){\line( 2,-3){600}}
\put(7501,-3661){\line( 0,-1){600}}
\put(7501,-4261){\line(-5,-3){1500}}
\put(6001,-2761){\line(-1,-1){900}}
\put(5101,-3661){\line( 3,-5){900}}
\put(6001,-4561){\line( 0,-1){600}}
\put(6001,-4561){\line( 5, 3){1500}}
\put(7501,-4261){\line(-2, 3){600}}
\put(6901,-2761){\line( 0,-1){600}}
\put(6001,-2761){\line( 0,-1){600}}
\put(6001,-3361){\line(-1,-1){900}}
\put(5101,-3661){\line( 0,-1){600}}
\put(5101,-4261){\line( 3,-1){900}}
\put(6001,-3361){\line( 1, 0){900}}
\put(6001,-2611){\makebox(0,0)[lb]{\smash{\SetFigFont{12}{14.4}{\rmdefault}{\mddefault}{\updefault}$a_2$}}}
\put(6901,-2611){\makebox(0,0)[lb]{\smash{\SetFigFont{12}{14.4}{\rmdefault}{\mddefault}{\updefault}$a_3$}}}
\put(7651,-3661){\makebox(0,0)[lb]{\smash{\SetFigFont{12}{14.4}{\rmdefault}{\mddefault}{\updefault}$a_4$}}}
\put(7651,-4336){\makebox(0,0)[lb]{\smash{\SetFigFont{12}{14.4}{\rmdefault}{\mddefault}{\updefault}$a_9$}}}
\put(6001,-5461){\makebox(0,0)[lb]{\smash{\SetFigFont{12}{14.4}{\rmdefault}{\mddefault}{\updefault}$a_{10}$}}}
\put(5926,-4336){\makebox(0,0)[lb]{\smash{\SetFigFont{12}{14.4}{\rmdefault}{\mddefault}{\updefault}$a_5$}}}
\put(4801,-4411){\makebox(0,0)[lb]{\smash{\SetFigFont{12}{14.4}{\rmdefault}{\mddefault}{\updefault}$a_6$}}}
\put(4801,-3661){\makebox(0,0)[lb]{\smash{\SetFigFont{12}{14.4}{\rmdefault}{\mddefault}{\updefault}$a_1$}}}
\put(6001,-3661){\makebox(0,0)[lb]{\smash{\SetFigFont{12}{14.4}{\rmdefault}{\mddefault}{\updefault}$a_7$}}}
\put(6676,-3661){\makebox(0,0)[lb]{\smash{\SetFigFont{12}{14.4}{\rmdefault}{\mddefault}{\updefault}$a_8$}}}
\end{picture}

\vspace*{0.5cm}

Figure 10
\end{center}

\vspace*{1cm}

{\bf Dynkin graphs}.
Let $D_n$ denote the tree with $n$ vertices as depicted in Figure 11.
For $n>4$ the automorphism group of $D_n$ is isomorphic to $Z_2$.
It has $n-1$ orbits on the set of vertices and $n^2-2n+2$ orbits on the set
of ordered pairs of vertices.

\noindent

\begin{center}
{\bf Table IV.} Results for Dynkin graphs

\vspace*{0.4cm}

\begin{tabular}{ccccc} \hline
   & number   & number    & \multicolumn{2}{c}{CPU time [seconds]} \\
 n & of cells & of colors & STABIL & STABCOL \\ \hline
6   & 5   & 26    &    0.03 &   0.01  \\
12  & 11  & 122   &    0.09 &   0.10  \\
18  & 17  & 290   &    0.14 &   0.36  \\
24  & 23  & 530   &    0.25 &   1.05  \\
30  & 29  & 842   &    0.43 &   2.28  \\
36  & 35  & 1226  &    0.79 &   5.23  \\
42  & 41  & 1682  &    1.34 &   9.94  \\
48  & 47  & 2210  &    1.95 &  15.84  \\
54  & 53  & 2810  &    2.78 &  23.82  \\
60  & 59  & 3482  &    4.04 &  33.08  \\
66  & 65  & 4226  &    5.30 &  44.21  \\
72  & 71  & 5042  &    7.55 &  57.97  \\ 
80  & 79  & 6242  &   10.66 &  82.19  \\
100 & 99  & 9802  &   26.72 & 171.07  \\
120 & 119 & 14162 &   83.12 & 303.23  \\
140 & 139 & 19322 &  258.99 & 488.47  \\
160 & 159 & 25282 &  729.08 &    --   \\
180 & 179 & 32042 & 1839.28 &    --   \\
\hline
\end{tabular}
\end{center}

\begin{center}
\setlength{\unitlength}{0.00083300in}%
\begingroup\makeatletter\ifx\SetFigFont\undefined%
\gdef\SetFigFont#1#2#3#4#5{%
  \reset@font\fontsize{#1}{#2pt}%
  \fontfamily{#3}\fontseries{#4}\fontshape{#5}%
  \selectfont}%
\fi\endgroup%
\begin{picture}(4666,1966)(5018,-4044)
\thicklines
\put(5101,-3061){\circle*{150}}
\put(6001,-3061){\circle*{150}}
\put(6901,-3061){\circle*{150}}
\put(7801,-3061){\circle*{150}}
\put(8701,-3061){\circle*{150}}
\put(9601,-2161){\circle*{150}}
\put(9601,-3961){\circle*{150}}
\put(5101,-3061){\line( 1, 0){1800}}
\multiput(6901,-3061)(30.00000,0.00000){31}{\makebox(6.6667,10.0000){\SetFigFont{10}{12}{\rmdefault}{\mddefault}{\updefault}.}}
\put(7801,-3061){\line( 1, 0){900}}
\put(9601,-2161){\line(-1,-1){900}}
\put(8701,-3061){\line( 1,-1){900}}
\end{picture}

\vspace*{0.5cm}

Figure 11 
\end{center}

The program STABIL has a rather long history. The first attempt of an
implementation was done by E.V.~Krukovskaya in PASCAL, see
\cite{KliK90}. A draft version of the present program was written in C by
I.V.~Chuvaeva and D.V.~Pasechnik at the N.~D.~Zelinski{\v i} Institute of
Organic Chemistry (Moscow) in 1990--1992, see \cite{ChuKP92}. Finally
this version was improved at the Technical University Munich in
1995. The improved version had static memory and, by this reason, was
available only to graphs with up to 200 vertices. This version was
carefully tested and the results of this testing are presented
above. In 1996 new improvements were done according to the suggestions
of Ch.~Pech (Dresden): dynamical memory management was created. Now the
current version, in principle, can handle graphs with an arbitrary
number of vertices. If for a given graph $\Gamma$ the number of
vertices is sufficiently small (that is if there will be enough memory
for saving all data structures) then we will get $W(\Gamma)$ as
output. Otherwise, the program will inform the user that the task cannot
be fulfilled completely. 

This last version of the program was used for other purely theoretical
purposes. Our experience shows that graphs with up to 500 vertices can
be successfully managed, however in some cases we were able to handle
even larger graphs. 

The codes of both programs STABIL and STABCOL and a read.me file which
describes how to use the programs are released under GPLv3
and can be downloaded from 
the homepage of the fourth author:  
\url{http://www.ntu.edu.sg/home/dima/software.htm}.

The users of the programs are requested to reference this report whenever
results  obtained with help  of STABIL or STABCOL are published.

\section{Discussion}\label{sect:disc}
The presented algorithm provides a very powerful and efficient tool to determine
equivalence of atoms and pairs of atoms in molecules. The equivalence classes are 
obtained by examining in a systematic way all configurations of three vertices in the 
underlying graph. The partitions of the vertices and edges are the finest which can 
be deduced using configurations of this size.

As already mentioned before, the standard partition of a graph not necessarily
coincides with its automorphism partition. Indeed, there exist graphs where the former 
partition is coarser than the latter. Methods which settle this shortcoming to
a certain extent are based on the following idea.
Classify vertices and edges by examining configurations which
consist not only of three but of a larger number of vertices.
This proceeding is generally called {\it deep stabilization\/} or
{\it stabilization of depth\/} $t$.

Roughly speaking, the situation is the following. Let $\Gamma$ be a graph with
$n$ vertices (possibly a directed multigraph) and $t$ a fixed integer, 
$2 \leq t \leq n$. All possible $n^t$ subgraphs of $\Gamma$ 
which are induced by the ordered $t$-tuples of vertices
are examined. We have to find all $\mu$ isomorphism types of these subgraphs.
To each pair $(u,v)$ of vertices a vector of length $\mu$ is associated, each
component of the vector being equal to the number of subgraphs of the corresponding
isomorphism type which contain the pair $(u,v)$. Any iteration of the stabilization
procedure of depth $t$ assigns two pairs $(u,v)$ and $(u',v')$ the same color if
and only if the vectors corresponding to these pairs are equal. 
It is clear that the computation of the total degree partition is nothing else than 
stabilization of depth $2$, Weisfeiler-Leman stabilization has depth $3$.
For depth at least 4 we obtain stronger algorithms, however at the price of a 
considerably higher complexity. One of the first attempts of a program
implementation of stabilization of depth $t\ge 4$ for purely chemical
goals and on a rather ``naive'' level was done in \cite{RueR91}.

It is demonstrated in the next paper \cite{FurKT} of this series that,
in contrast to first expectations, stabilization of depth $t$ with some 
$t \geq 4$ is also not sufficient to rigorously settle the automorphism 
partitioning problem. It turns out (see \cite{Fur87}, \cite{CaiFI92}) that
for any fixed value of $t$ there exist graphs with the
property that the standard partition of depth $t$ does not coincide with the 
automorphism partition.

In future work we intend to develop an implementation of the WL-stabilization
which eventually is even faster than the implementation STABIL presented in
this paper. A very promising approach is to perform in an alternative way
stabilization steps of depth 2 and 3. Given a colored complete directed graph
$\Delta=(\Omega,R)$, we start with stabilization of depth 2, i.e.\ we compute
the total degree partition of $\Delta$ (let $\Delta_k$ be the graph consisting
of the arcs of color $k$; then the total degree partition of $\Delta$ is the
coarsest partition of $V$ such that any two vertices belonging to the same cell
of the partition have the same valencies with respect to any other cell in any
graph $\Delta_k$). The arcs of a given color are recolored according to the
colors of their end vertices such that arcs between different colored pairs of
vertices are assigned different colors. In the next step the coloring of the
arcs and vertices is refined analogously as in STABIL by considering all
triangles of the graph. However, in order to decrease the effort, this is not
done iteratively, but only once. After that we again compute the total degree
partition, perform one stabilization step of depth 3, etc. The algorithm stops
if the coloring of $\Delta$ is stable.


As it was mentioned before, we can only be sure to get the automorphism partition
of a graph $\Gamma$  by means of WL-stabilization if it is known in advance that
the algebra $W(\Gamma)$ is Schurian.

In general, we can only suggest to proceed in the following way:\\
find the automorphism group $G=\aut(\Gamma)$ of the graph $\Gamma$;\\
describe the set of 2--orbits (or only 1--orbits) of the action of $G$
on the vertex set of $\Gamma$.

This problem, in principle, may be solved using e.g. the computer
package COCO (I.A.~Farad\v{z}ev, M.H.~Klin), the UNIX implementation by
A.E.~Brouwer.

The preliminary versions of COCO are described in \cite{FarK91} and
\cite{FarKM94}. With the use of COCO one may handle graphs with a few
thousands of vertices. The current version of COCO \cite{COCO} 
 is oriented towards purely
mathematical goals, namely for the investigation of graphs having a
sufficiently large prescribed subgroup of the automorphism group. In
such a case the input graph is described by a set of arcs as a
union of suitable 2--orbits of a prescribed permutation group. For
purely chemical purposes such a mode of input is certainly inconvenient.
Hopefully, in the future a more suitable interface for chemists 
will be created.

This technical report is considered by the authors as a preliminary
version of a future regular publication. We will be very grateful to
everybody who supplies us with remarks, comments, criticism or
improvements.

Putting the original version of \cite{BabCKP97} onto {\sf arXiV},
we just updated few references, while not trying to reflect
the ongoing progress during the last decade. The reader may benefit
from recent papers \cite{EvdPT00,CohKM08,EvdP09}, where some
new relevant ideas and methods are discussed. The authors still
hope to produce a throughout revision of the present text.

Unless more convenient, M.~Klin should be regarded as the corresponding
author.
   
\parskip0.0cm
\section*{Acknowledgments}
We want to express our deep gratitude to the following colleagues:

G.M.~Adel'son--Velski\v{\i} for helpful discussions concerning the complexity of
the WL--stabilization;

S.~Baumann and M.~L\"udecke for their help in preparing the
manuscript and testing the program;

I.A.~Farad\v{z}ev for very important consultations related to dealing
with hashing problems (during the creation of the first draft of the
program);

Ch.~Pech for all suggested improvements, in particular for assistance in implementing
a dynamical memory management;

Ch. and G.~R\"ucker for the fruitful discussions and remarks which
reflect real requests of chemists to the use of algebraic
combinatorics;

G.~Tinhofer for permanent interest and support of this project;

N.S.~Zefirov for creating a lucky opportunity for mathematicians to
work with chemists and to realize their practical interests.

Last but not the least, G.~Cooperman for pointing out a number of typos
in \cite{BabCKP97}.

\end{document}